\documentclass[onefignum,onetabnum]{siamart190516}

\usepackage{amsmath,amsfonts,amssymb,mathtools,nicefrac,mymath}
\usepackage{graphicx,rotating}
\usepackage[caption=false]{subfig}
\usepackage[algo2e,ruled,linesnumbered]{algorithm2e}
\usepackage{comment}
\usepackage{soul,xcolor,cancel}
\usepackage{braket}
\usepackage{cite}
\usepackage{enumitem}

\newcommand\cpp{C\nolinebreak[4]\hspace{-.05em}\raisebox{.4ex}{\relsize{-3}{\textbf{++}}}}

\definecolor{myblue}{rgb}{0,0.4470,0.7410}
\definecolor{myred}{rgb}{0.8500,0.3250,0.0980}
\definecolor{myorange}{rgb}{0.9290,0.6940,0.1250}
\definecolor{mypurple}{rgb}{0.4940,0.1840,0.5560}
\definecolor{mygreen}{rgb}{0.4660,0.6740,0.1880}
\definecolor{mylightblue}{rgb}{0.3010,0.7450,0.9330}
\definecolor{mydarkred}{rgb}{0.6350,0.0780,0.1840}

  \usepackage{tikz}
  \usepackage{pgfplots}
  \usetikzlibrary{calc}
  \usetikzlibrary{positioning}
  \usetikzlibrary{shapes.multipart}
  \usetikzlibrary{pgfplots.groupplots}
  \usepackage{mycdgates}
  \pgfplotsset{
    compat=newest,
    table/header=false,
    tick label style={font=\scriptsize},
    label style={font=\scriptsize},
    legend style={font=\scriptsize},
    legend cell align=left,
    colormap={parula}{
      rgb255=(53,42,135)
      rgb255=(15,92,221)
      rgb255=(18,125,216)
      rgb255=(7,156,207)
      rgb255=(21,177,180)
      rgb255=(89,189,140)
      rgb255=(165,190,107)
      rgb255=(225,185,82)
      rgb255=(252,206,46)
      rgb255=(249,251,14)
    }
  }
  \pgfplotsset{
    myColOne/.style={myblue},
    myColTwo/.style={myred},
    myColThr/.style={myorange},
    myColFou/.style={mypurple},
    myColFiv/.style={mygreen},
    myColSix/.style={mylightblue},
    myColSev/.style={mydarkred}
  }
  \pgfplotsset{
    myStyOne/.style={myColOne,thick,mark=+},
    myStyTwo/.style={myColTwo,thick,mark=+},
    myStyThr/.style={myColThr,thick,mark=+},
    myStyFou/.style={myColFou,thick,mark=+},
    myStyFiv/.style={myColFiv,thick,mark=+},
    myStySix/.style={myColSix,thick,mark=+},
    myStySev/.style={myColSev,thick,mark=+},
    myStyRes/.style={black,densely dotted},
  }
  \pgfplotsset{
    myStyOn/.style={myblue,     thick,mark=asterisk},
    myStyTw/.style={myred,      thick,mark=asterisk},
    myStyTh/.style={myorange,   thick,mark=asterisk},
    myStyFo/.style={mypurple,   thick,mark=asterisk},
    myStyFi/.style={mygreen,    thick,mark=asterisk},
    myStySi/.style={mylightblue,thick,mark=asterisk},
    myStySe/.style={mydarkred,  thick,mark=asterisk},
  }
  \pgfplotsset{
    every axis/.append style={
      label style={font=\footnotesize},
    }
  }
  \newcommand{\figname}[1]{\relax}
  \newcommand{\fignames}[1]{\relax}
  \newcommand{\datfile}[1]{fig/#1.dat}

\usepackage{qcircuit}
\usepackage{environ}
\newcommand{\myqctmp}[2][0.5]{\Qcircuit @C=#2em @R=#1em @!R}
\NewEnviron{myqcircuit}[1][0.5]{\vcenter{\myqctmp[#1]{0.8} {\BODY}}}
\NewEnviron{myqcircuit*}[2]{\vcenter{\myqctmp[#1]{#2} {\BODY}}}
\newcommand{\qdots}{\push{~\cdots~}\qw}

\usepackage[binary-units]{siunitx}
\usepackage{tabularx,multirow,booktabs}
\newcolumntype{L}{>{\raggedright\arraybackslash}X}
\newcolumntype{C}{>{\centering\arraybackslash}X}
\newcolumntype{R}{>{\raggedleft\arraybackslash}X}

\crefname{subsection}{section}{sections}
\Crefname{subsection}{Section}{Sections}
\newsiamthm{assumption}{Assumption}
\crefname{assumption}{Assumption}{Assumptions}
\newsiamthm{conjecture}{Conjecture}
\crefname{conjecture}{Conjecture}{Conjectures}
\newsiamremark{example}{Example}
\crefname{example}{Example}{Examples}
\newsiamremark{remark}{Remark}
\crefname{remark}{Remark}{Remarks}
\crefname{algocf}{Algorithm}{Algorithms}

\if@review
  \usepackage{etoolbox}
  \newcommand*\linenomathpatchAMS[1]{%
    \expandafter\pretocmd\csname #1\endcsname {\linenomathAMS}{}{}%
    \expandafter\pretocmd\csname #1*\endcsname{\linenomathAMS}{}{}%
    \expandafter\apptocmd\csname end#1\endcsname {\endlinenomath}{}{}%
    \expandafter\apptocmd\csname end#1*\endcsname{\endlinenomath}{}{}%
  }
  \expandafter\ifx\linenomath\linenomathWithnumbers
    \let\linenomathAMS\linenomathWithnumbers
    \patchcmd\linenomathAMS{\advance\postdisplaypenalty\linenopenalty}{}{}{}
  \else
    \let\linenomathAMS\linenomathNonumbers
  \fi
  \linenomathpatchAMS{gather}
  \linenomathpatchAMS{multline}
  \linenomathpatchAMS{align}
  \linenomathpatchAMS{alignat}
  \linenomathpatchAMS{flalign}
\fi

\DeclareMathOperator{\SU2}{SU(2)}
\DeclareMathOperator{\TFXY}{TFXY}

\newcommand{\TheTitle}%
{An Algebraic Quantum Circuit Compression Algorithm for Hamiltonian Simulation}
\newcommand{\TheTitleArticle}%
{An Algebraic Quantum Circuit Compression Algorithm for Hamiltonian Simulation}
\newcommand{\TheAuthors}%
{D.~Camps, E.~K\"okc\"u, et al.}

\headers{An algebraic circuit compression algorithm}{\TheAuthors}

\title{{\TheTitleArticle}\thanks{Submitted to the editors \today.
\funding{DC and RVB are supported by the Laboratory Directed Research and Development Program of Lawrence Berkeley National Laboratory under U.S. Department of Energy Contract No. DE-AC02-05CH11231.
LB, and WAdJ were supported by the U.S. Department of Energy (DOE) under Contract No.  DE-AC02-05CH11231,  through the Office of Advanced Scientific Computing  Research  Accelerated  Research  for  Quantum Computing Program. EK, and AFK were supported by the Department of Energy, Office of Basic Energy Sciences, Division of Materials Sciences and Engineering under Grant No. DE-SC0019469.
}}}

\author{
    Daan~Camps%
    \thanks{Computational Research Division,
            Lawrence Berkeley National Laboratory,
            Berkeley, CA 94720, United States.
            (\email{dcamps@lbl.gov},
             \email{lbassman@lbl.gov},
             \email{wadejong@lbl.gov},
             \email{rvanbeeumen@lbl.gov}).}
    \and         
  	Efekan~K\"{o}kc\"{u}%
    \thanks{Department of Physics, 
          North Carolina State University, 
          Raleigh, NC 27695, United States.
          (\email{ekokcu@ncsu.edu},
           \email{akemper@ncsu.edu}).}           
  \and
    Lindsay~Bassman%
    \footnotemark[2]
  \and
    Wibe~A.~de~Jong%
    \footnotemark[2]
  \and
    Alexander~F.~Kemper%
    \footnotemark[3]  
  \and
    Roel~Van~Beeumen%
    \footnotemark[2]
}

\ifpdf
\hypersetup{
  pdftitle={\TheTitle},
  pdfauthor={\TheAuthors}
}
\fi

\begin{document}
\maketitle


\begin{abstract}
Quantum computing is a promising technology that harnesses the peculiarities
of quantum mechanics to deliver computational speedups for some problems that
are intractable to solve on a classical computer.
Current generation noisy intermediate-scale quantum (NISQ) computers are severely
limited in terms of chip size and error rates.
Shallow quantum circuits with uncomplicated topologies are essential for successful
applications in the NISQ era.
Based on matrix analysis, we derive localized circuit transformations to efficiently 
compress quantum circuits for simulation of certain
spin Hamiltonians known as free fermions.
The depth of the compressed circuits is independent of simulation time and grows 
linearly with the number of spins.
The proposed numerical circuit compression algorithm behaves backward stable and
scales cubically in the number of spins enabling circuit synthesis beyond $\bigO(10^3)$ spins.
The resulting quantum circuits have a simple nearest-neighbor topology, which makes them
ideally suited for NISQ devices.
\end{abstract}




\section{Introduction}
\label{sec:intro}

The field of quantum computing~\cite{NC2010} is rapidly evolving.
The current generation of quantum hardware is known as 
\emph{Noisy Intermediate-Scale Quantum} (NISQ) computers~\cite{Preskill2018}
and can perform specialized computational tasks that become rapidly intractable
for a classical computer~\cite{Google}.
A computational task for a quantum computer, or \emph{quantum program}, is typically
expressed as a \emph{quantum circuit} that consists of a sequence
of unitary transformations~\cite{NC2010}. 
Each of these unitary transformations typically acts on just one or two \emph{qubits} 
of the quantum computer and they are often referred to as \emph{quantum gates}.
The technological limitations in NISQ hardware impose substantial constraints 
both on the number of qubits and on the number of
unitary operations, also known as \emph{circuit depth}, that can be performed.
Noise introduced by two qubit gates eventually reduces the fidelity of the quantum state until
a useful signal can no longer be measured.
Shallow and simple quantum circuit structures are thus crucial for successful applications
in the NISQ era.

Quantum circuit \emph{compilation} or \emph{synthesis}~\cite{Barenco1995a,Shende} is
the problem of computing a circuit representation into two qubit operations 
for a target unitary matrix.
General purpose synthesis algorithms based on well-known matrix decompositions have 
been proposed in the literature, for example, based on 
Givens~\cite{Vartiainen2004} or Householder~\cite{Ivanov2006} QR factorization 
of the target unitary.
A more efficient algebraic synthesis algorithm in terms of circuit 
complexity is known as the \emph{quantum Shannon decomposition}~\cite{Shende}
and is based on a hierarchical cosine-sine decomposition (CSD)~\cite{Sutton2009,Sutton2012}
of the unitary matrix.
While these methods work for every unitary matrix, they have two major disadvantages.
Firstly, they require an exorbitant amount of classical resources.
The dimension of the unitary matrix for $N$ qubits to be decomposed is $2^N \times 2^N$.
Storing this matrix on a classical computer
rapidly becomes intractable, let alone computing a decomposition of cubic complexity
in the matrix dimension such as a QR factorization or CSD.
Secondly, the circuits that are derived from these synthesis algorithms contain,
in general, exponentially many gates in terms of the number of qubits.
For many applications of practical interest, more efficient circuits can be obtained
with optimization methods~\cite{qfast} or by exploiting certain structures in the unitary~\cite{Camps2020}.
This approach still suffers from the first issue as the full unitary has to be formed.

In this paper we propose an application-specific circuit compression and synthesis
algorithm that overcomes both challenges. We never form the $2^N \times 2^N$
unitary and the compression algorithm has a cubic complexity in $N$
which is an exponential improvement compared to a cubic dependence on $2^N$~\cite{Vartiainen2004,Ivanov2006,Shende}.
Furthermore, the compressed circuits have a simple nearest-neighbor topology, a circuit
depth of $\bigO(N)$, and $\bigO(N^2)$ quantum gates.
This makes them ideally suited for the NISQ era and in particular for hardware
based on superconducting qubits.
The application that our synthesis algorithm is designed for is known as \emph{Hamiltonian
simulation}~\cite{Lloyd96} which involves the evolution of a quantum state of the
system under the time-dependent Schr\"odinger equation.
This problem is ubiquitous in quantum chemistry~\cite{Bassman2021a,Bauer2020}
and physics, for example in adiabatic ground state preparation~\cite{Barends2016}.
We show that quantum circuits for the time evolution of certain spin models, 
known in physics as free fermionizable or integrable models, are efficiently compressible.
Our analysis leads to an algebraic circuit compression algorithm that
behaves as a backward stable algorithm in practice.
MATLAB and \cpp\ implementations of our algorithms are publicly available
as part of the \emph{fast free fermion compiler} (\texttt{F3C}) \cite{f3c,f3cpp} at 
\url{https://github.com/QuantumComputingLab}.
\texttt{F3C} is build based on the \texttt{QCLAB} toolbox \cite{qclab,qclabpp} for creating and
representing quantum circuits.

A related algorithm based on a Givens QR factorization of an $N \times N$ matrix
formed by the quadratic Hamiltonian was proposed in \cite{Kivlichan2018} to generate
Slater determinants.
This method was further developed to prepare a Hartree-Fock wave function \cite{Arute2020}
and generic fermionic Gaussian states \cite{Jiang2018}.
It assumes that the full circuit maps to a free fermionic system,
while our localized operations can still be used for circuits that are partially
comprised of specific quantum gates.

This paper is accompanied by a dual paper targeted
at the physics community~\cite{TrotterCompression} that analyzes the properties 
of the Hamiltonians from studying the Hamiltonian algebra.
While this paper focuses on the matrix structures and the efficient
and accurate numerical computation of the compression, \cite{TrotterCompression}
focuses on the implications for the physics community and showcases the
results of an adiabatic state preparation experiment performed on quantum
hardware that is only feasible due to the compressed quantum circuits.

The remainder of our paper is organized as follows. 
\Cref{sec:problem} provides a more detailed introduction to the problem of
Hamiltonian simulation, reviews the concept of operator splitting methods
to solve this problem, introduces the specific spin Hamiltonians for which
our compression algorithm works, and relates the current paper to earlier
work.
\Cref{sec:props} provides an overview of useful elementary results on Pauli rotation
matrices and parameterizations of $\SU2$ that we will use
for the remainder of the analysis.
\Cref{sec:ising} shows that compressing
quantum circuits for the simulation of classical Ising models to depth $\bigO(1)$
immediately follows from \Cref{sec:props}.
We present our circuit compression algorithms for simulation circuits
that are comprised of gates that allow for a \emph{fusion} and \emph{turnover} operation
in \Cref{sec:algorithms}.
\Cref{sec:KXY} demonstrates based on the results from \Cref{sec:props} that Kitaev
chains and XY Hamiltonians satisfy these criteria and can be efficiently compressed.
\Cref{sec:transverse} shows the same for transverse-field XY Hamiltonians and the special
case of transverse-field Ising models.
Details of the implementation and considerations on numerical stability are
provided in \Cref{sec:implementation}. 
\Cref{sec:ex} provides numerical examples that demonstrate the speed and accuracy
of our method. We conclude in \Cref{sec:concl}.

\section{Problem statement and preliminary results} 
\label{sec:problem}

In this section we review the problem statement and preliminary
results about Hamiltonian simulation on quantum computers.

\subsection{Hamiltonian simulation}
\label{ssec:hamsim}

Simulating a quantum system of $N$ \emph{spins} or \emph{qubits}
involves the evolution of the quantum state of the system 
under the Schr\"odinger equation,
\begin{equation*}
\frac{\partial}{\partial t} \psi(t) = - \I H(t) \psi(t),
\end{equation*}
and is fully determined by $H(t) \in \C^{2^N \times 2^N}$, 
the time-dependent Hamiltonian of the system, 
and the initial state of the system, $\psi(t_0) \in \C^{2^N}$.
The Hamiltonian is a time-dependent Hermitian operator of exponential
dimension in the system size and the quantum state $\psi(t)$ is
a vector of unit norm.

Simulating from initial time $t_0$ to final time $t_1$ 
is achieved by the time-evolution operator
\begin{equation}
    U(t_1,t_0)=\mathcal{T} \exp\left(-\I\int_{t_0}^{t_1}H(t)dt\right),
    \label{eq:timeevol}
\end{equation}
where $\mathcal{T} \exp$ is the \emph{time-ordered matrix exponential}.
The final state at time $t_1$ becomes
$\psi(t_1) = U(t_1,t_0) \psi(t_0)$.
For a time-independent Hamiltonian the closed-form solution is
$U(t_1,t_0) = \exp\left( -\I (t_1 - t_0) H\right)$.
Time evolution is a hard problem to solve on a classical computer due to the
exponential dimensionality of the state space and the time-dependence 
of the Hamiltonian.
In \emph{digital quantum simulation}, \eqref{eq:timeevol} is evaluated on 
a quantum computer which naturally operates in a state space of exponential
dimension.
Multiple quantum algorithms have been proposed with (near) optimal asymptotic
scaling~\cite{Low2017,Low2019,Gilyen2018b,Berry2015a,Berry2015c,Kalev2020,Haah2021}.
All of these algorithms rely on more complicated circuit structures that are
not well-suited for the constraints imposed by NISQ computers where the circuit
depth is limited.
Quantum circuits derived from
\emph{operator splitting}~\cite{McLachlan2002,Thalhammer2012}
methods, also known as
\emph{Trotter product formulas}~\cite{trotter,suzuki,Hatano2005,Childs2021}
in the physics community, often lead to simple circuit structures but with a
circuit depth that usually depends linearly on simulation time.

\subsection{Operator splitting methods}
\label{sec:prodform}

We rely on two approximations in order to implement \eqref{eq:timeevol} 
on a quantum computer using an operator splitting method.
First, we discretize in time by approximating $H(t)$ by a piecewise
constant function with $n_t$ time-steps of length $\Delta t$ that discretize
the interval $[t_0, t_1)$~\cite{Poulin2011}:
\begin{align*}
H(t) & \approx H(t_\tau) =: H_{\tau}, & 
0 < \tau & \leq n_t, & 
t_\tau & = t_0 + \tau \Delta t, &
t & \in [t_{\tau - 1}, t_\tau). 
\end{align*}

Second, we approximate the matrix exponential of $H_{\tau}$ by 
products of matrix exponentials that are easier to implement
on a quantum computer.
The simplest case is a first-order product formula which decomposes
the Hamiltonian operator in a sum of two terms $H = A + B$.
The approximate time-evolution operator for time-step $\Delta t$,
$U(\Delta t) = \exp(-\I A \Delta t) \exp(-\I B \Delta t)$, satisfies~\cite{Childs2021}:
\begin{equation*}
\| U(\Delta t) - \exp(-\I H \Delta t) \| \leq \frac{\Delta t^2}{2} \| [A, B] \|,
\end{equation*}
where $[A, B] := AB - BA$.
This result can be bootstrapped to show that for $H = \sum_{\ell} H_{\ell}$ and
$U(\Delta t) = \prod_{\ell} \exp(-\I H_\ell \Delta t)$, we have that
\begin{equation}
\label{eq:Trottererr}
\| U(\Delta t) - \exp(-\I H \Delta t) \| \leq \frac{\Delta t^2}{2} \sum_{i>j} \| [H_i, H_j] \|.
\end{equation}
Without loss of generality,
we will only use first-order Trotter decompositions throughout this paper.

Combining the discretization in time ($H(t) \approx H_\tau$) 
and the Trotter decomposition of the
Hamiltonian, $H_\tau = \sum_{\ell} H_{\ell,\tau}$, we
have the following approximation to the time-evolution operator
\begin{align}
U(n_t \Delta t) & = \prod_{k=0}^{n_t -1} U_{n_t - k}(\Delta t), &
U_{\tau}(\Delta t) & = \prod_{\ell} \exp(-\I H_{\ell,\tau}\Delta t), 
\label{eq:dtimeevol}
\end{align}
where the index $\tau$ is a time-ordered multiplication over $n_t$ discretized
time-steps $\Delta t$ and $\ell$ multiplies over the terms in the Trotter
decomposition.
Quantum circuits based on this formula naturally become a concatenation
of blocks that implement the individual time-steps and their depth grows
linearly with $n_t$:
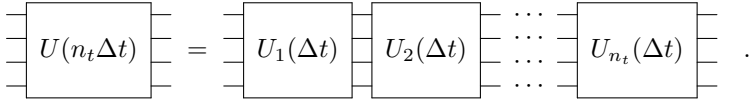
\begin{equation*}
{\begin{myqcircuit*}{0}{0.75}
  & \multigate3{U(n_t \Delta t)} & \qw \\
  & \ghost{U(n_t \Delta t)} & \qw \\
  & \ghost{U(n_t \Delta t)} & \qw \\
  & \ghost{U(n_t \Delta t)} & \qw
\end{myqcircuit*}}
\ = \
{\begin{myqcircuit*}{0}{0.75}
  & \multigate3{U_1(\Delta t)} & \multigate3{U_2(\Delta t)} & \qdots & \multigate3{U_{n_t}(\Delta t)} & \qw \\
  & \ghost{U_1(\Delta t)} & \ghost{U_2(\Delta t)} & \qdots & \ghost{U_{n_t}(\Delta t)} & \qw \\
  & \ghost{U_1(\Delta t)} & \ghost{U_2(\Delta t)} & \qdots & \ghost{U_{n_t}(\Delta t)} & \qw \\
  & \ghost{U_1(\Delta t)} & \ghost{U_2(\Delta t)} & \qdots & \ghost{U_{n_t}(\Delta t)} & \qw
\end{myqcircuit*}} \quad.
\end{equation*}

Note that in quantum circuit diagrams, time flows from left to right, which means
that the order of operations is reversed compared to \eqref{eq:dtimeevol}.

\subsection{Spin Hamiltonians}
\label{ssec:spinham}

Our compression results only hold for certain time-dependent, \emph{ordered} or \emph{disordered}
Hamiltonians that model chains of
spin-$\nicefrac{1}{2}$ particles with a nearest-neighbor coupling and
external magnetic field.
These Hamiltonians are typically expressed in terms of the Pauli spin-$\nicefrac{1}{2}$ matrices,
\begin{align*}
\sigma^{x} & = \begin{bmatrix} 0 & 1 \\ 1 & 0 \end{bmatrix}, &
\sigma^{y} & = \begin{bmatrix} 0 & -\I \\ \I & \phantom{-}0 \end{bmatrix}, &
\sigma^{z} & = \begin{bmatrix} 1 & \phantom{-}0 \\ 0 & -1 \end{bmatrix},
\end{align*}
which are generators for the group of $2 \times 2$ unitary matrices with unit determinant, also known as $\SU2$.
We often write $\sigma^{\alpha}$ where $\alpha \in \lbrace x, y, z \rbrace$
as many of our results are independent of the type of Pauli matrix.
A basis for the Hilbert space of composite quantum systems is constructed through the tensor
product of the state spaces of the individual systems.
To this end it is useful to introduce an abbreviated notation
for a $\sigma^{\alpha}$ matrix that acts on the $i$th spin in a chain of
$N$ spins:
\begin{equation}
\sigma^{\alpha}_i := \underbrace{\eye \otimes \cdots \otimes \eye}_{i-1} \, \otimes \, \sigma^{\alpha} \otimes \underbrace{\eye \otimes \cdots \otimes \eye}_{N-i},
\end{equation}
where $\eye$ is the $2 \times 2$ identity matrix.

The most complicated class of Hamiltonians that we consider is known as 
the \emph{time-dependent, disordered transverse field XY} (TFXY) model. 
The Hamiltonian is given by
\begin{align}
H(t) = 
& \underbrace{\sum_{i=1}^{N-1} J^{x}_i(t) \, \sigma_{i}^{x} \sigma_{i+1}^{x} + J^{y}_i(t)  \,\sigma_{i}^{y} \sigma_{i+1}^{y}}_{\text{Coupling}}
\ + \ \underbrace{\sum_{i=1}^{N} h^{z}_i(t) \, \sigma_{i}^z}_{\text{External Field}}.
\label{eq:TFXY}
\end{align}
The other two permutations of $x$, $y$, and $z$ result in TFXZ and
TFYZ Hamiltonians for which our circuit compression method also works.

The parameters $J^{x}$, $J^{y}$, and $h^{z}$ in the Hamiltonian \eqref{eq:TFXY} depend 
both on time $t$ and on the index $i$ of the spin-$\nicefrac{1}{2}$ particle in the chain. 
The dependence on the index $i$ means that the Hamiltonian is disordered. 
If the parameters are independent of $i$, the Hamiltonian is called ordered.
We discuss TFXY Hamiltonians in detail in \Cref{sec:transverse}.
Other models that are compressible are classical Ising models, Kitaev chains,
XY models, and tranverse-field Ising models (TFIM).
All of these are subclasses of TFXY Hamiltonians obtained by restricting
some parameters of the full TFXY model.

\subsection{Related work}
\label{ssec:related}

Besides~\cite{TrotterCompression}, this paper is closely related to two earlier
papers~\cite{Bassman2021,Kokcu2021} written by some of us.
The results in our current paper were first conjectured in~\cite{Bassman2021}.
There, the fixed-depth circuit property was identified by using QFAST~\cite{qfast}, an optimization-based numerical
circuit compiler. This compilation method becomes challenging for problems larger than
$7$ qubits because of the exponential dimensionality of the state space and this
prompted us to analyze the problem in more detail. 
Our analysis presented in the current paper resulted in a constructive proof of the fixed-depth
property (\Cref{sec:algorithms}) for all cases conjectured in~\cite{Bassman2021}
and a scalable and accurate circuit compression algorithm that easily handles systems with
$\bigO(10^3)$ qubits.
Furthermore, we improve the circuit depth of the classical Ising model to $\bigO(1)$
compared to~\cite{Bassman2021}.
In~\cite{Kokcu2021} the existence of fixed depth circuits for Hamiltonian 
simulation is proven through a Cartan decomposition of the Lie algebra
generated by the Hamiltonian.
The advantage of our current approach over~\cite{Kokcu2021} is that by compressing
an easy to generate Trotter circuit to shallow depth, we avoid having to optimize
the whole circuit at once, which again only scales up to $10$ qubits~\cite{Kokcu2021}.
Furthermore, our results are derived from matrix analysis and our algorithms
are exact up to machine precision.


\section{Elementary properties and results}
\label{sec:props}

We give an overview of all properties that we use later in the analysis of 
the circuit compression algorithm.
We start with the following well-known commutation relations between Pauli
matrices:
\begin{align}
\label{eq:paulicom}
\left[\sigma^\alpha, \sigma^\beta\right] & = 2 \I \varepsilon_{\alpha\beta\gamma} \sigma^\gamma, 
& \lbrace \alpha, \beta, \gamma \rbrace \subseteq \lbrace x, y, z \rbrace,
\end{align}
with $\varepsilon_{\alpha\beta\gamma}$ the Levi-Cevita tensor.
Kronecker products of two Pauli matrices do commute:
\begin{align}
\label{eq:2paulicom}
\left[\sigma^\alpha \otimes \sigma^\alpha, \sigma^\beta \otimes \sigma^\beta\right] & = 0,
& \alpha, \beta \in \lbrace x, y, z \rbrace.
\end{align}

\begin{definition}\label{def:rotations}
For $\alpha \in \lbrace x, y, z \rbrace$,
we define a single-spin Pauli-$\alpha$ rotation over an angle $\theta$ as
\begin{equation}
R^{\alpha}(\theta) := \exp(-\I \, \sigma^\alpha \, \theta/2 ) =
\begin{tikzpicture}[on grid, auto, baseline=-0.5ex]
\fqwire{w1}
\rextgate[myblue]{gate1}{w11}{$\alpha$}
\end{tikzpicture}.
\end{equation}
If the Pauli-$\alpha$ rotation acts on the $i$th spin, we denote it as $R^{\alpha}_{i}(\theta)$.
Similarly, a two-spin Pauli-$\alpha$ rotation over an angle $\theta$ is defined as
\begin{equation}
R^{\alpha\alpha}(\theta) := \exp(-\I \, \sigma^\alpha \otimes \sigma^\alpha \, \theta/2 )
=
\begin{tikzpicture}[on grid, auto,baseline=-0.5ex]
\fnode{node}{}
\arqwires[2*\HorDist][\VerDist][\VerDist/2 and 0mm]{w}{node}{2}
\broagate[myblue][\VerDist/2 and \HorDist]{gate}{w11}{$\alpha$}
\end{tikzpicture},
\end{equation}
and denoted as $R^{\alpha\alpha}_{i}(\theta)$ if it acts on the nearest-neighbor spins $i, i+1$.
\end{definition}

We introduced our diagrammatic notation for single- and two-spin Pauli-rotations
in \Cref{def:rotations}.
This will be our quantum circuit representation for these unitary matrices.
The vertical direction
indicates the spins on which the unitary operations are performed.
The matrix representation of the single- and two-spin Pauli rotations are given by:
\begin{align*}
R^x(\theta) & = \begin{bmatrix} c & -\I s\\ -\I s & c \end{bmatrix}, &
R^{xx}(\theta) & = \begin{bmatrix} c &  &  & -\I s\\  &c & -\I s &  \\  & -\I s & c & \\ -\I s &  &  &c \end{bmatrix},\\
R^y(\theta) & = \begin{bmatrix} c & -s\\ s & c \end{bmatrix}, &
R^{yy}(\theta) & = \begin{bmatrix} c & & & \I s\\  & c & - \I s & \\  & -\I s & c &  \\ \I s & & & c \end{bmatrix},\\
R^z(\theta) & = \begin{bmatrix} e^{-\I \theta/2} & \\ & e^{\I \theta/2}\end{bmatrix}, &
R^{zz}(\theta) & =  \begin{bmatrix} e^{-\I \theta/2} & & & \\  & e^{\I \theta/2} &  & \\ & & e^{\I \theta/2} & \\ & & & e^{-\I \theta/2}\end{bmatrix},
\end{align*}
where $c = \cos(\nicefrac{\theta}{2})$ and $s = \sin(\nicefrac{\theta}{2})$.

A property that we often use implicitly is the mixed product property
of the Kronecker product and the observation that the identity commutes with every
matrix:
\begin{align}
\label{eq:commute1}
(A \otimes I)(I \otimes B) & = (I \otimes B)(A \otimes I), &
\begin{tikzpicture}[on grid, auto, baseline=-\VerDist/2]
\fqwires[3*\HorDist]{w1}{2}
\rextgate[black][2*\HorDist]{g1}{w111}{$A$}
\rextgate[black]{g2}{w121}{$B$}
\arnode[\VerDist/2 and \HorDist/2]{eq1}{w122}{$=$}
\rqwires[3*\HorDist]{w2}{w112}{2}
\rextgate[black]{g3}{w211}{$A$}
\rextgate[black][2*\HorDist]{g4}{w221}{$B$}
\end{tikzpicture}
\end{align}
Here we have a first illustration of the reversed order of operations
in matrix notation compared to the schematic notation.

The following three lemmas list useful properties of the Pauli rotation matrices 
that are directly verified from \Cref{def:rotations} in combination with the commutation
relations \eqref{eq:paulicom,eq:2paulicom,eq:commute1}.
The first lemma provides some useful commutation relations for Pauli rotation
matrices.

\begin{lemma}\label{lem:pcom}
For $\alpha, \beta \in \lbrace x, y, z \rbrace$, $i \in \lbrace 1, \dots, N-1 \rbrace$,
the following commutation relations hold for the Pauli rotations:
\begin{enumerate}[wide, label=(\roman*),ref=(\roman*)]
\item{\label{itm:pcomaa} Two-spin rotations of the same type on overlapping spins:
\begin{align*}
R^{\alpha\alpha}_{i+1}(\theta_2) \, R^{\alpha\alpha}_{i}(\theta_1) & = R^{\alpha\alpha}_{i}(\theta_1) \, R^{\alpha\alpha}_{i+1}(\theta_2), &
\begin{tikzpicture}[on grid, auto, baseline=-\VerDist]
\fqwires[3*\HorDist]{w3}{3}
\broagate[myblue][\VerDist/2 and \HorDist]{g5}{w311}{$\alpha$}
\broagate[myblue][\VerDist and \HorDist]{g6}{g5}{$\alpha$}
\rnode[\HorDist/4]{eq2}{w322}{$=$}
\rqwires[3*\HorDist][\VerDist][\HorDist/2]{w4}{w312}{3}
\aroagate[myblue][\VerDist/2 and \HorDist]{g7}{w431}{$\alpha$}
\aroagate[myblue][\VerDist and \HorDist]{g8}{g7}{$\alpha$}
\lnode[\HorDist/2]{l1}{w311}{\footnotesize$i$}
\lnode[\HorDist/2]{l2}{w321}{\footnotesize$i{+}1$}
\lnode[\HorDist/2]{l3}{w331}{\footnotesize$i{+}2$}
\bnode[9*\HorDist/4]{p1}{g5}{\footnotesize$\theta_1$}
\bnode[5*\HorDist/4]{p2}{g6}{\footnotesize$\theta_2$}
\bnode[5*\HorDist/4]{p3}{g7}{\footnotesize$\theta_2$}
\bnode[9*\HorDist/4]{p4}{g8}{\footnotesize$\theta_1$}
\end{tikzpicture}
\end{align*}
}
\item{\label{itm:pcomaabb} Two-spin rotations of different type on the same spins:
\begin{align*}
R^{\beta\beta}_{i}(\theta_2) R^{\alpha\alpha}_{i}(\theta_1) & = R^{\alpha\alpha}_{i}(\theta_1) \, R^{\beta\beta}_{i}(\theta_2), &
\begin{tikzpicture}[on grid, auto, baseline=-\VerDist/2]
\fqwires[3*\HorDist]{w1}{2}
\broagate[myblue][\VerDist/2 and \HorDist]{g1}{w111}{$\alpha$}
\roagate[myred]{g2}{g1}{$\beta$}
\arnode[\VerDist/2 and \HorDist/4]{eq1}{w122}{$=$}
\rqwires[3*\HorDist][\VerDist][0.5*\HorDist]{w2}{w112}{2}
\broagate[myred][\VerDist/2 and \HorDist]{g3}{w211}{$\beta$}
\roagate[myblue]{g4}{g3}{$\alpha$}
\lnode[\HorDist/2]{l1}{w111}{\footnotesize$i$}
\lnode[\HorDist/2]{l2}{w121}{\footnotesize$i{+}1$}
\bnode[5*\HorDist/4]{p1}{g1}{\footnotesize$\theta_1$}
\bnode[5*\HorDist/4]{p2}{g2}{\footnotesize$\theta_2$}
\bnode[5*\HorDist/4]{p3}{g3}{\footnotesize$\theta_2$}
\bnode[5*\HorDist/4]{p4}{g4}{\footnotesize$\theta_1$}
\end{tikzpicture}
\end{align*}
}
\item{\label{itm:pcomaaa} Single- and two-spin rotations of the same type on the same spins:
\begin{align*}
R^{\alpha}_{i}(\theta_3) R^{\alpha}_{i+1}(\theta_2) \, R^{\alpha\alpha}_{i}(\theta_1) & = R^{\alpha\alpha}_{i}(\theta_1) \, R^{\alpha}_{i+1}(\theta_2) R^{\alpha}_{i}(\theta_3), &
\begin{tikzpicture}[on grid, auto, baseline=-\VerDist/2]
\fqwires[3*\HorDist]{w1}{2}
\broagate[myblue][\VerDist/2 and \HorDist]{g4}{w111}{$\alpha$}
\lextgate[myblue]{g5}{w112}{$\alpha$}
\lextgate[myblue]{g6}{w122}{$\alpha$}
\arnode[\VerDist/2 and \HorDist/4]{eq1}{w122}{$=$}
\rqwires[3*\HorDist][\VerDist][0.5*\HorDist]{w2}{w112}{2}
\rextgate[myblue]{g1}{w211}{$\alpha$}
\rextgate[myblue]{g2}{w221}{$\alpha$}
\broagate[myblue][\VerDist/2 and \HorDist]{g3}{g1}{$\alpha$}
\lnode[\HorDist/2]{l1}{w111}{\footnotesize$i$}
\lnode[\HorDist/2]{l2}{w121}{\footnotesize$i{+}1$}
\anode[3*\HorDist/4]{p1}{g1}{\footnotesize$\theta_3$}
\bnode[3*\HorDist/4]{p2}{g2}{\footnotesize$\theta_2$}
\bnode[5*\HorDist/4]{p3}{g3}{\footnotesize$\theta_1$}
\anode[3*\HorDist/4]{p5}{g5}{\footnotesize$\theta_3$}
\bnode[3*\HorDist/4]{p6}{g6}{\footnotesize$\theta_2$}
\bnode[5*\HorDist/4]{p4}{g4}{\footnotesize$\theta_1$}
\end{tikzpicture}
\end{align*}
}
\end{enumerate}
\end{lemma}

The next lemma shows that rotations of the same type acting
on the same spins can be easily fused together.
\begin{lemma}\label{lem:pfus}
Let $\alpha \in \lbrace x, y, z \rbrace$. 
\begin{enumerate}[wide,label=(\roman*),ref=(\roman*)]
\item{\label{itm:pfusa} For $i \in \lbrace 1, \dots, N \rbrace$, single-spin Pauli-$\alpha$ rotations acting on the same spins 
can be \emph{fused} or multiplied together:
\begin{align*}
R^{\alpha}_{i}(\theta_2) \, R^{\alpha}_{i}(\theta_1) & = R^{\alpha}_{i}(\theta_1 + \theta_2), &
\begin{tikzpicture}[on grid, auto, baseline=-\VerDist/4]
\fqwire[3*\HorDist]{w1}
\rextgate[myblue]{g1}{w11}{$\alpha$}
\rextgate[myblue]{g2}{g1}{$\alpha$}
\rnode[\HorDist/2]{eq1}{w12}{$=$}
\rqwire{w2}{w12}
\rextgate[myblue]{g3}{w21}{$\alpha$}
\lnode[\HorDist/2]{l1}{w11}{\footnotesize$i$}
\bnode[3*\HorDist/4]{p1}{g1}{\footnotesize$\theta_1$}
\bnode[3*\HorDist/4]{p2}{g2}{\footnotesize$\theta_2$}
\bnode[3*\HorDist/4]{p3}{g3}{\footnotesize$\theta_1 + \theta_2$}
\end{tikzpicture}
\end{align*}
}
\item{\label{itm:fusaa} For $i \in \lbrace 1, \dots, N-1 \rbrace$, two-spin Pauli-$\alpha$ rotations acting on the same spins 
can be \emph{fused} or multiplied together:
\begin{align*}R^{\alpha\alpha}_{i}(\theta_2) \, R^{\alpha\alpha}_{i}(\theta_1) & = R^{\alpha\alpha}_{i}(\theta_1 + \theta_2),&
\begin{tikzpicture}[on grid, auto, baseline=-\VerDist/2]
\fqwires[3*\HorDist]{w1}{2}
\broagate[myblue][\VerDist/2 and \HorDist]{g4}{w11}{$\alpha$}
\roagate[myblue]{g5}{g4}{$\alpha$}
\brnode[\VerDist/2 and \HorDist/2]{eq2}{w112}{$=$}
\rqwires{w2}{w112}{2}
\broagate[myblue][\VerDist/2 and \HorDist]{g6}{w211}{$\alpha$}
\lnode[\HorDist/2]{l1}{w111}{\footnotesize$i$}
\lnode[\HorDist/2]{l2}{w121}{\footnotesize$i{+}1$}
\bnode[5*\HorDist/4]{p1}{g4}{\footnotesize$\theta_1$}
\bnode[5*\HorDist/4]{p2}{g5}{\footnotesize$\theta_2$}
\bnode[5*\HorDist/4]{p3}{g6}{\footnotesize$\theta_1 + \theta_2$}
\end{tikzpicture}
\end{align*}
}
\end{enumerate}
\end{lemma}

The following result directly follows from the commutativity \eqref{eq:2paulicom}.
\begin{lemma}
\label{lem:2axes}
For $\alpha, \beta \in \lbrace x, y, z \rbrace$, $\alpha \neq \beta$, $0 \leq \theta_\alpha, \theta_\beta < 4\pi$, we have that
\begin{align*}
\exp\left(-\I \left(\sigma^\alpha \otimes \sigma^\alpha \theta_\alpha/2 + \sigma^\beta \otimes \sigma^\beta \theta_\beta/2 \right)\right)
& = \exp(-\I \sigma^\alpha \otimes \sigma^\alpha \theta_\alpha/2) \exp(-\I \sigma^\beta \otimes \sigma^\beta \theta_\beta/2),\\
R^{\alpha\beta}(\theta_\alpha,\theta_\beta)
& = R^{\alpha\alpha}(\theta_\alpha) R^{\beta\beta}(\theta_\beta),
\end{align*}
or as a circuit diagram:
\begin{equation*}
\begin{tikzpicture}[on grid, auto]
\fqwires{w1}{2}
\brtagate[\VerDist/2 and \HorDist]{myblue,myred}{gate1}{w111}{$\alpha$}{$\beta$}
\arnode[\VerDist/2 and \HorDist/4]{eq1}{w122}{$=$}
\rqwires[3*\HorDist][\VerDist][0.5*\HorDist]{w2}{w112}{2}
\broagate[myblue][\VerDist/2 and \HorDist]{gate2}{w211}{$\alpha$}
\roagate[myred]{gate3}{gate2}{$\beta$}
\arnode[\VerDist/2 and \HorDist/4]{eq1}{w222}{$=$}
\rqwires[3*\HorDist][\VerDist][0.5*\HorDist]{w3}{w212}{2}
\broagate[myred][\VerDist/2 and \HorDist]{gate4}{w311}{$\beta$}
\roagate[myblue]{gate5}{gate4}{$\alpha$}
\end{tikzpicture}
\end{equation*}
\end{lemma}

\subsection{Euler decompositions of SU(2)}

The group $\SU2$ is given by
\begin{equation}
\SU2 = \left\{ \begin{bmatrix} \alpha & -\bar{\beta}\\ \beta & \bar{\alpha} \end{bmatrix}: \alpha, \beta \in \C, |\alpha|^2 + |\beta|^2 = 1 \right\},
\label{eq:SU2}
\end{equation}
and it is well-known (see for example \cite[Theorem 4.1]{NC2010}) that any element of $\SU2$
can be parametrized by three \emph{Euler angles}. We will refer to this as an
\emph{Euler decomposition} of $\SU2$.

\begin{lemma}
Let $\alpha, \beta \in \lbrace x, y, z \rbrace$, $\alpha \neq \beta$. Every matrix $U \in \SU2$
can be represented as:
\begin{align}
U & = R^{\alpha}(\theta_1) \, R^{\beta}(\theta_2) \, R^{\alpha}(\theta_3), &
\begin{tikzpicture}[on grid, auto, baseline=-\VerDist/4]
\fqwire{w1}
\rextgate[black]{g1}{w11}{$U$}
\rnode[\HorDist/2]{eq1}{w12}{$=$}
\rqwire[4*\HorDist]{w2}{w12}
\rextgate[myblue]{g2}{w21}{$\alpha$}
\rextgate[myred]{g3}{g2}{$\beta$}
\rextgate[myblue]{g4}{g3}{$\alpha$}
\bnode[3*\HorDist/4]{p1}{g2}{\footnotesize$\theta_3$}
\bnode[3*\HorDist/4]{p2}{g3}{\footnotesize$\theta_2$}
\bnode[3*\HorDist/4]{p3}{g4}{\footnotesize$\theta_1$}
\end{tikzpicture}.
\end{align}
The decomposition is unique, except for a set of measure 0, if the angles are restricted to $0 \leq \theta_1 < 2\pi$,
$0 \leq \theta_2 \leq \pi$ and $0 \leq \theta_3 < 4\pi$.
\end{lemma}

\begin{proof}
We give a proof for $\alpha = z$, $\beta = y$, other
cases follow from a similar argument.
Direct computation yields:
\begin{align}
R^{z}(\theta_1) \, R^{y}(\theta_2) \, R^{z}(\theta_3) & =
\begin{bmatrix} e^{-\I \theta_1/2} & \\ & e^{\I \theta_1/2} \end{bmatrix}
\begin{bmatrix} \cos(\nicefrac{\theta_2}{2}) & -\sin(\nicefrac{\theta_2}{2}) \\ \sin(\nicefrac{\theta_2}{2}) &\phantom{-}\cos(\nicefrac{\theta_2}{2}) \end{bmatrix}
\begin{bmatrix} e^{-\I \theta_3/2} & \\ & e^{\I \theta_3/2} \end{bmatrix}, \nonumber\\
& =
\begin{bmatrix}
\cos(\nicefrac{\theta_2}{2}) e^{-\I (\theta_1 + \theta_3)/2} & -\sin(\nicefrac{\theta_2}{2}) e^{-\I (\theta_1 - \theta_3)/2} \\ \sin(\nicefrac{\theta_2}{2}) e^{\phantom{-}\I (\theta_1 - \theta_3)/2} & \phantom{-}\cos(\nicefrac{\theta_2}{2}) e^{\phantom{-}\I (\theta_1 + \theta_3)/2}
\end{bmatrix}.
\label{eq:zyz}
\end{align}
It is clear that this parametrizes $\SU2$ and that the Euler angles are unique unless $\theta_2 = 0, \pi$.
\end{proof}

We can \emph{turn over} an Euler decomposition of $\SU2$ to its dual decomposition as shown in
the following result.

\begin{lemma}\label{lem:turnoverSU2}
Let $\alpha, \beta \in \lbrace x, y, z \rbrace$, $\alpha \neq \beta$.
For every set of Euler angles given by $\theta_1$, $\theta_2$, $\theta_3$ there exists a set of dual 
Euler angles given by $\theta_a$, $\theta_b$, $\theta_c$ such that,
\begin{align*}
R^{\alpha}(\theta_1) \, R^{\beta}(\theta_2) \, R^{\alpha}(\theta_3) & = R^{\beta}(\theta_a) \, R^{\alpha}(\theta_b) \, R^{\beta}(\theta_c), &
\begin{tikzpicture}[on grid, auto, baseline=-\VerDist/4]
\fqwire[4*\HorDist]{w1}
\rextgate[myblue]{g1}{w11}{$\alpha$}
\rextgate[myred]{g2}{g1}{$\beta$}
\rextgate[myblue]{g3}{g2}{$\alpha$}
\rnode[\HorDist/4]{eq1}{w12}{$=$}
\rqwire[4*\HorDist][\HorDist/2]{w2}{w12}
\rextgate[myred]{g4}{w21}{$\beta$}
\rextgate[myblue]{g5}{g4}{$\alpha$}
\rextgate[myred]{g6}{g5}{$\beta$}
\bnode[3*\HorDist/4]{p1}{g1}{\footnotesize$\theta_3$}
\bnode[3*\HorDist/4]{p2}{g2}{\footnotesize$\theta_2$}
\bnode[3*\HorDist/4]{p3}{g3}{\footnotesize$\theta_1$}
\bnode[3*\HorDist/4]{p1}{g4}{\footnotesize$\theta_c$}
\bnode[3*\HorDist/4]{p2}{g5}{\footnotesize$\theta_b$}
\bnode[3*\HorDist/4]{p3}{g6}{\footnotesize$\theta_a$}
\end{tikzpicture}.
\label{eq:Eulerto}
\end{align*}
The relation between both sets of Euler angles is given by:
\begin{equation}
\label{eq:toSU21}
\begin{aligned}
\tan\left(\frac{\theta_a + \theta_c}{2}\right) & =
\phantom{-}\tan\left(\frac{\theta_2}{2}\right) \, \frac{\cos((\theta_1 - \theta_3)/2)}{\cos((\theta_1 + \theta_3)/2)}, \\
\tan\left(\frac{\theta_a - \theta_c}{2}\right) & =
-\tan\left(\frac{\theta_2}{2}\right) \, \frac{\sin((\theta_1 - \theta_3)/2)}{\sin((\theta_1 + \theta_3)/2)},
\end{aligned}
\end{equation}
and
\begin{equation}
\label{eq:toSU22}
\begin{aligned}
\tan\left(\frac{\theta_1 + \theta_3}{2}\right) & =
\phantom{-}\tan\left(\frac{\theta_b}{2}\right) \, \frac{\cos((\theta_a - \theta_c)/2)}{\cos((\theta_a + \theta_c)/2)}, \\
\tan\left(\frac{\theta_1 - \theta_3}{2}\right) & =
-\tan\left(\frac{\theta_b}{2}\right) \, \frac{\sin((\theta_a - \theta_c)/2)}{\sin((\theta_a + \theta_c)/2)}.
\end{aligned}
\end{equation}
\end{lemma}

\begin{proof}
We again only give the proof for $\alpha = z$, $\beta = y$, but the result holds in general.
In this case, the left-hand side of \Cref{eq:Eulerto} is given by \eqref{eq:zyz}.
The right-hand side is equal to:
\begin{align*}
& \begin{bmatrix} \cos(\nicefrac{\theta_a}{2}) & -\sin(\nicefrac{\theta_a}{2}) \\ \sin(\nicefrac{\theta_a}{2}) &\phantom{-}\cos(\nicefrac{\theta_a}{2}) \end{bmatrix}
\begin{bmatrix} e^{-\I \theta_b/2} & \\ & e^{\I \theta_b/2} \end{bmatrix} 
\begin{bmatrix} \cos(\nicefrac{\theta_c}{2}) & -\sin(\nicefrac{\theta_c}{2}) \\ \sin(\nicefrac{\theta_c}{2}) &\phantom{-}\cos(\nicefrac{\theta_c}{2}) \end{bmatrix}.
\end{align*}
As this is an $\SU2$ matrix, it is determined by its first column. 
After performing the matrix products and using some elementary trigonometry relations, we find the first column to be:
\begin{equation}
\begin{aligned}
& (1,1): & \cos(\nicefrac{\theta_b}{2}) \cos((\theta_a + \theta_c)/2) - \I \sin(\nicefrac{\theta_b}{2}) \cos((\theta_a - \theta_c)/2), \\
& (2,1): & \cos(\nicefrac{\theta_b}{2}) \sin((\theta_a + \theta_c)/2) - \I \sin(\nicefrac{\theta_b}{2}) \sin((\theta_a - \theta_c)/2).
\end{aligned}
\label{eq:yzy}
\end{equation}
The lemma directly follows from setting the first column of this matrix to the first column of \eqref{eq:zyz}.
\end{proof}

The following two lemmas are closely related to \Cref{lem:turnoverSU2}.

\begin{lemma}\label{lem:mato}
Let $\alpha, \beta \in \lbrace x, y, z \rbrace$, $\alpha \neq \beta$, $i \in \lbrace 1, \dots, N-1 \rbrace$.
For every set of Euler angles $\theta_1, \theta_2,\theta_3$,
there exists a set of dual Euler angles $\theta_a, \theta_b, \theta_c$ such that
\begin{equation*}
R^{\alpha\alpha}_{i}(\theta_1) \, R^{\beta\beta}_{i+1}(\theta_2) \, R^{\alpha\alpha}_{i}(\theta_3)  = R^{\beta\beta}_{i+1}(\theta_a) \, R^{\alpha\alpha}_{i}(\theta_b) \, R^{\beta\beta}_{i+1}(\theta_c).
\end{equation*}
In terms of a diagram this relation is given by:
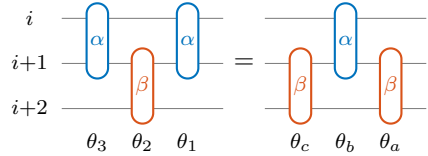
\begin{equation*}
\begin{tikzpicture}[on grid, auto]
\fqwires[4*\HorDist]{w1}{3}
\broagate[myblue][\VerDist/2 and \HorDist]{gate3}{w111}{$\alpha$}
\broagate[myred]{gate2}{gate3}{$\beta$}
\aroagate[myblue]{gate1}{gate2}{$\alpha$}
\rnode[\HorDist/4]{eq1}{w122}{$=$}
\rqwires[4*\HorDist][\VerDist][0.5*\HorDist]{w2}{w112}{3}
\broagate[myred][\VerDist/2 and \HorDist]{gateC}{w221}{$\beta$}
\aroagate[myblue]{gateB}{gateC}{$\alpha$}
\broagate[myred]{gateA}{gateB}{$\beta$}
\lnode[\HorDist/2]{l1}{w111}{\footnotesize$i$}
\lnode[\HorDist/2]{l2}{w121}{\footnotesize$i{+}1$}
\lnode[\HorDist/2]{l2}{w131}{\footnotesize$i{+}2$}
\bnode[9*\HorDist/4]{p1}{gate3}{\footnotesize$\theta_3$}
\bnode[5*\HorDist/4]{p2}{gate2}{\footnotesize$\theta_2$}
\bnode[9*\HorDist/4]{p3}{gate1}{\footnotesize$\theta_1$}
\bnode[5*\HorDist/4]{p1}{gateC}{\footnotesize$\theta_c$}
\bnode[9*\HorDist/4]{p2}{gateB}{\footnotesize$\theta_b$}
\bnode[5*\HorDist/4]{p3}{gateA}{\footnotesize$\theta_a$}
\end{tikzpicture}
\end{equation*}
The relation between both sets of Euler angles is given by \eqref{eq:toSU21,eq:toSU22}.
\end{lemma}

\begin{lemma}\label{lem:TFIMto}
Let $\alpha, \beta \in \lbrace x, y, z \rbrace$, $\alpha \neq \beta$, $i \in \lbrace 1, \dots, N-1 \rbrace$.
For every set of Euler angles $\theta_1, \theta_2,\theta_3$,
there exists a set of dual Euler angles $\theta_a, \theta_b, \theta_c$ such that
\begin{align*}
\label{eq:TFIMto}
R^{\alpha\alpha}_{i}(\theta_1) \, R^{\beta}_{i}(\theta_2) \, R^{\alpha\alpha}_{i}(\theta_3)  & = R^{\beta}_{i}(\theta_a) \, R^{\alpha\alpha}_{i}(\theta_b) \, R^{\beta}_{i}(\theta_c), \\
R^{\alpha\alpha}_{i}(\theta_1) \, R^{\beta}_{i+1}(\theta_2) \, R^{\alpha\alpha}_{i}(\theta_3)  & = R^{\beta}_{i+1}(\theta_a) \, R^{\alpha\alpha}_{i}(\theta_b) \, R^{\beta}_{i}(\theta_c),
\end{align*}
or as a circuit diagram:
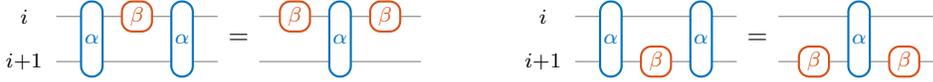
\begin{equation*}
\begin{tikzpicture}[on grid, auto]
\fqwires[4*\HorDist]{w1}{2}
\broagate[myblue][\VerDist/2 and \HorDist]{g4}{w111}{$\alpha$}
\rextgate[myred][2*\HorDist]{g5}{w111}{$\beta$}
\bloagate[myblue][\VerDist/2 and \HorDist]{g4}{w112}{$\alpha$}
\brnode[\VerDist/2 and \HorDist/4]{eq1}{w112}{$=$}
\rqwires[4*\HorDist][\VerDist][0.5*\HorDist]{w2}{w112}{2}
\rextgate[myred]{g1}{w211}{$\beta$}
\broagate[myblue][\VerDist/2 and 2*\HorDist]{g2}{w211}{$\alpha$}
\rextgate[myred][2*\HorDist]{g3}{g1}{$\beta$}
\lnode[\HorDist/2]{l1}{w111}{\footnotesize$i$}
\lnode[\HorDist/2]{l2}{w121}{\footnotesize$i{+}1$}
\rqwires[4*\HorDist][\VerDist][3*\HorDist]{w3}{w212}{2}
\broagate[myblue][\VerDist/2 and \HorDist]{g10}{w311}{$\alpha$}
\rextgate[myred][2*\HorDist]{g11}{w321}{$\beta$}
\bloagate[myblue][\VerDist/2 and \HorDist]{g12}{w312}{$\alpha$}
\brnode[\VerDist/2 and \HorDist/4]{eq2}{w312}{$=$}
\rqwires[4*\HorDist][\VerDist][0.5*\HorDist]{w4}{w312}{2}
\rextgate[myred]{g7}{w421}{$\beta$}
\broagate[myblue][\VerDist/2 and 2*\HorDist]{g8}{w411}{$\alpha$}
\rextgate[myred][2*\HorDist]{g9}{g7}{$\beta$}
\lnode[\HorDist/2]{l1}{w311}{\footnotesize$i$}
\lnode[\HorDist/2]{l2}{w321}{\footnotesize$i{+}1$}
\end{tikzpicture}
\end{equation*}
The relation between both sets of Euler angles is given by \eqref{eq:toSU21,eq:toSU22}.
\end{lemma}

\Cref{lem:mato,lem:TFIMto} follow from the observation that the matrices involved
have the same group structure as $\SU2$. We refer the interested reader to \cite{TrotterCompression}
for further details.


\section{Classical Ising model}
\label{sec:ising}

In this section we show, based on the results from \Cref{sec:props},
that the approximate time-evolution operator from
\eqref{eq:dtimeevol} can be implemented in a circuit of depth $\bigO(1)$
for Hamiltonians that are known as (classical) Ising models.
The Ising model is a classical Hamiltonian because all the terms in
its expansion commute. This allows for the $\bigO(1)$ depth.
The Hamiltonians in \Cref{sec:KXY,sec:transverse} consist of terms
that do not commute and are truly quantum. They will
require deeper circuits that are more challenging to compute.

The Hamiltonian for the Ising model is given by,
\begin{align}
\label{eq:Hising}
H(t) & = \sum_{i=1}^{N-1} J^{\alpha}_i(t) \, \sigma^{\alpha}_i \sigma^{\alpha}_{i+1} + \sum_{i=1}^{N} h^{\alpha}_i(t) \, \sigma^{\alpha}_i, &\alpha\in \lbrace x, y, z \rbrace.
\end{align}
We will use the shorthand notation $H^{\alpha\alpha + \alpha}(t)$ for this Hamiltonian thereby
referring to its nonzero terms.
If we decompose this Hamiltonian in its two-spin and single-spin interaction, $H^{\alpha\alpha + \alpha}(t) = H^{\alpha\alpha}(t) + H^{\alpha}(t)$, a single time-step in the discretized time-evolution operator \eqref{eq:dtimeevol} becomes,
\begin{align*}
U_{\tau}(\Delta t) & = \exp(-\I H^{\alpha\alpha}_{\tau}\Delta t) \, \exp(-\I H^{\alpha}_{\tau}\Delta t), \\
& = \exp\left(-\I  \sum_{i=1}^{N-1} J^{\alpha}_i(t_{\tau}) \, \sigma^{\alpha}_i \sigma^{\alpha}_{i+1} \, \Delta t \right) \, \exp\left(-\I \sum_{i=1}^{N} h^{\alpha}_i(t_{\tau}) \, \sigma^{\alpha}_i \, \Delta t\right), \\
& = \prod_{i=1}^{N-1} R^{\alpha\alpha}_{i}(2 J^{\alpha}_i(t_\tau) \, \Delta t) \ \prod_{i=1}^{N} R^{\alpha}_{i}(2 h^{\alpha}_i(t_\tau) \, \Delta t).
\end{align*}

Since all the terms in $H^{\alpha\alpha + \alpha}(t)$ commute according to \eqref{eq:paulicom,eq:commute1},
we did not introduce a Trotter error according to \eqref{eq:Trottererr}, except for the discretization in time.
Using the commutativity of two-spin Pauli rotations, \Cref{lem:pcom}\ref{itm:pcomaa}, it follows
that we can rearrange the ascending cascades of two-spin Pauli rotations that make up a single time-step
into an even-odd ordering:
\begin{align*}
\prod_{i=1}^{N-1} R^{\alpha\alpha}_{i}(2 J^{\alpha}_i(t_\tau)  \Delta t) = 
\prod_{\text{odd } i} R^{\alpha\alpha}_{i}(2 J^{\alpha}_i(t_\tau)  \Delta t)
\prod_{\text{even } i} R^{\alpha\alpha}_{i}(2 J^{\alpha}_i(t_\tau)  \Delta t)
\end{align*}
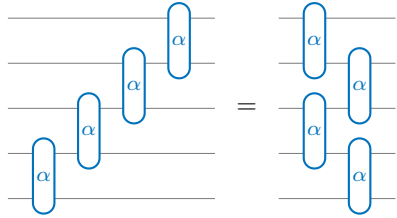
\begin{center}
\begin{tikzpicture}[on grid,auto]
\fqwires[5*\HorDist]{w1}{5}
\anode[\VerDist/2]{n1}{w151}{}
\roaslantedlayer[myblue][\VerDist and \HorDist][\HorDist]{l1}{n1}{$\alpha$}{4}
\rqwires[3*\HorDist]{w2}{w112}{5}
\broalayer[myblue][\DblVerDist][\VerDist/2 and \HorDist]{l12}{w211}{$\alpha$}{2}
\broalayer[myblue]{l13}{l121}{$\alpha$}{2}
\rnode[0.5*\HorDist]{eq}{w132}{$=$}
\end{tikzpicture}
\end{center}

\begin{figure}[hbtp]
\centering
\begin{tikzpicture}[on grid,auto]
\fqwires[6.5*\HorDist]{w1}{5}
\rqwires[2*\HorDist][\VerDist][-0.5*\HorDist][\WireColor,dashed]{w2}{w112}{5}
\rqwires[3.5*\HorDist][\VerDist][-0.5*\HorDist]{w3}{w212}{5}
\rqwires[4*\HorDist]{w4}{w312}{5}
\rextlayer[myblue]{l1}{w111}{$\alpha$}{5}
\broalayer[myblue][\DblVerDist][\VerDist/2 and \HorDist]{l2}{l11}{$\alpha$}{2}
\broalayer[myblue]{l3}{l21}{$\alpha$}{2}
\rextlayer[myblue][\VerDist][3*\HorDist]{l4}{l11}{$\alpha$}{5}
\broalayer[myblue][\DblVerDist][\VerDist/2 and \HorDist]{l5}{l41}{$\alpha$}{2}
\broalayer[myblue]{l6}{l51}{$\alpha$}{2}
\rextlayer[myblue][\VerDist][4*\HorDist]{l7}{l41}{$\alpha$}{5}
\broalayer[myblue][\DblVerDist][\VerDist/2 and \HorDist]{l8}{l71}{$\alpha$}{2}
\broalayer[myblue]{l9}{l81}{$\alpha$}{2}
\rextlayer[myblue]{l10}{w411}{$\alpha$}{5}
\broalayer[myblue][\DblVerDist][\VerDist/2 and \HorDist]{l11}{l101}{$\alpha$}{2}
\broalayer[myblue]{l12}{l111}{$\alpha$}{2}
\rnode[0.5*\HorDist]{eq}{w332}{$=$}
\bnode{lbl1}{l15}{\footnotesize$\theta_{:,1}$}
\bnode{lbl2}{l45}{\footnotesize$\theta_{:,2}$}
\bnode{lbl3}{l75}{\footnotesize$\theta_{:,n_t}$}
\blnode[\VerDist and 0.05*\HorDist]{lbl4}{l105}{\footnotesize$\sum_i \theta_{:,i}$}
\blnode[1.5*\VerDist and \HorDist/2]{lbl5}{l32}{\footnotesize$\phi_{:,1}$}
\blnode[1.5*\VerDist and \HorDist/2]{lbl6}{l62}{\footnotesize$\phi_{:,2}$}
\blnode[1.5*\VerDist and \HorDist/2]{lbl7}{l92}{\footnotesize$\phi_{:,n_t}$}
\blnode[1.5*\VerDist and 0.45*\HorDist]{lbl8}{l122}{\footnotesize$\sum_i \phi_{:,i}$}
\end{tikzpicture}
\vspace{-5pt}
\label{fig:pfcompression}
\caption{Compression to a circuit of depth $\bigO(1)$ for time-evolution of an Ising model Hamiltonian as in \eqref{eq:Hising}.}
\end{figure}
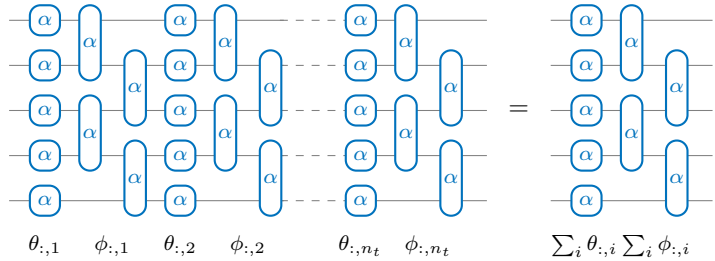

The time-evolution operator \eqref{eq:dtimeevol} for the Ising becomes a horizontal concatenation
of these layers interleaved with layers of single-spin rotations.
According to \Cref{lem:pcom} all gates in the circuit commute with each other and according
to \Cref{lem:pfus}, gates of the same type that are acting on the same spin(s) can be fused together. 
This circuit compression method for Ising models is illustrated in \Cref{fig:pfcompression}.

The parameters of the compressed $\bigO(1)$ circuits are computed by 
straightforward numerical integration of the parameters of the 
time-dependent Ising Hamiltonians.
This requires $\bigO(n_t N)$ operations for disordered Ising Hamiltonians
and only $\bigO(n_t)$ for ordered Ising Hamiltonians.
In case the Ising Hamiltonian is time-independent, the complexity can be 
further reduced to $\bigO(\log_2(n_t) N)$ and $\bigO(\log_2(n_t))$ for respectively 
disordered and ordered Ising Hamiltonians by applying the merging algorithm recursively.


\section{Constant-depth circuits with fusion and turnover operations}
\label{sec:algorithms}

In this section we present a second algorithm for computing
constant-depth circuits for time-evolution by means of 
compressing longer circuits.
The difference with the circuits presented in \Cref{sec:ising}
is the type of operations that we can perform on
the gates in the circuit.
We constructively show that quantum circuits
comprised of non-commuting two-spin gates which allow for a fusion and turnover operation
can efficiently be compressed to a circuit with depth $\bigO(N)$.
To show this, we use circuit transformations that are equivalent to 
transformations used in core chasing algorithms for eigenvalue 
problems~\cite{Aurentz2018,Vandebril2011,Vandebril2012}.

The turnover operation acts locally on a pattern of three two-spin
gates and changes a $\VEE$-shaped pattern to a $\HAT$-shaped pattern 
or vice versa:
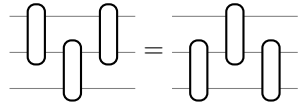
\begin{equation*}
\begin{tikzpicture}[on grid, auto]
\renewcommand{\IntHeight}{8mm}
\renewcommand{\HorDist}{4.8mm}
\renewcommand{\VerDist}{4.8mm}
\renewcommand{\DblVerDist}{2*\VerDist} 
\renewcommand{\DblHorDist}{2*\HorDist}
\tikzstyle{oneaxis}=[draw,rectangle,thick,rounded corners = 1mm,minimum height=\IntHeight,adjust width=\IntAspect,inner sep=0,fill=white]
\fqwires[4*\HorDist]{w1}{3}
\broagate[black][\VerDist/2 and \HorDist]{gate3}{w111}{}
\broagate[black]{gate2}{gate3}{}
\aroagate[black]{gate1}{gate2}{}
\rnode[\HorDist/4]{eq1}{w122}{$=$}
\rqwires[4*\HorDist][\VerDist][0.5*\HorDist]{w2}{w112}{3}
\broagate[black][\VerDist/2 and \HorDist]{gateC}{w221}{}
\aroagate[black]{gateB}{gateC}{}
\broagate[black]{gateA}{gateB}{}
\end{tikzpicture}
\end{equation*}
\Cref{lem:mato} was a first example of a turnover operation for 
two-spin Pauli rotations.
A useful operation that we can do with the turnover operation is 
illustrated in \Cref{fig:gtl1,fig:gtl2,fig:gtl3,fig:gtl4,fig:gtl5}.
It pulls a \emph{free gate} through an \emph{ascending cascade} of gates
which moves the incoming gate on position down.
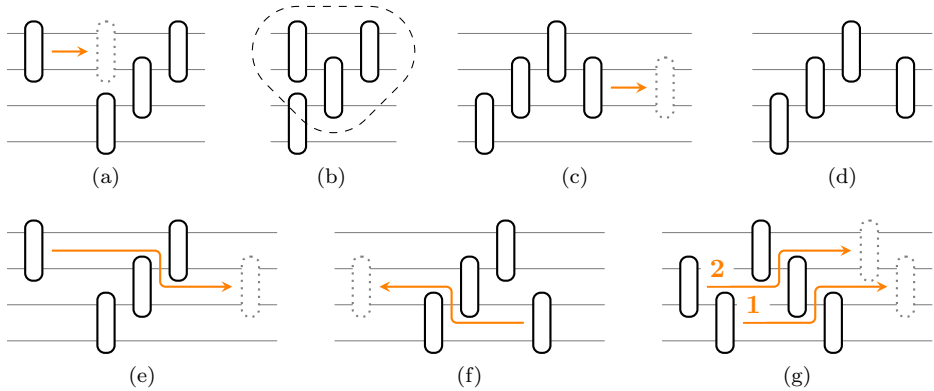
\begin{figure}[hbtp]
\centering
\subfloat[]{%
\label{fig:gtl1}
\begin{tikzpicture}[on grid,auto]
\renewcommand{\IntHeight}{8mm}
\renewcommand{\HorDist}{4.8mm}
\renewcommand{\VerDist}{4.8mm}
\renewcommand{\DblVerDist}{2*\VerDist} 
\renewcommand{\DblHorDist}{2*\HorDist}
\tikzstyle{oneaxis}=[draw,rectangle,thick,rounded corners = 1mm,minimum height=\IntHeight,adjust width=\IntAspect,inner sep=0,fill=white]
\fqwires[6*\HorDist]{w1}{4}
\anode[\VerDist/2]{node1}{w141}{}
\roaslantedlayer[black][\VerDist and \HorDist][3*\HorDist]{l1}{node1}{}{3}
\broagate[black][\VerDist/2 and \HorDist]{g1}{w111}{}
\aoagate[gray,dotted]{ga}{l11}{}
\coordinate[right= \HorDist/2 of g1] (a1);
\coordinate[right= \HorDist of a1] (a2);
\draw [myarrow] (a1) -- (a2);
\end{tikzpicture}
}\hfill%
\subfloat[]{%
\label{fig:gtl2}
\begin{tikzpicture}[on grid,auto]
\renewcommand{\IntHeight}{8mm}
\renewcommand{\HorDist}{4.8mm}
\renewcommand{\VerDist}{4.8mm}
\renewcommand{\DblVerDist}{2*\VerDist} 
\renewcommand{\DblHorDist}{2*\HorDist}
\tikzstyle{oneaxis}=[draw,rectangle,thick,rounded corners = 1mm,minimum height=\IntHeight,adjust width=\IntAspect,inner sep=0,fill=white]
\fqwires[4*\HorDist]{w2}{4}
\anode[\VerDist/2]{node2}{w241}{}
\roaslantedlayer[black][\VerDist and \HorDist][\HorDist]{l2}{node2}{}{3}
\broagate[black][\VerDist/2 and \HorDist]{g2}{w211}{}
\draw[dashed] \convexpath{l22,g2,l23}{6mm};
\end{tikzpicture}
}\hfill%
\subfloat[]{%
\label{fig:gtl3}
\begin{tikzpicture}[on grid,auto]
\renewcommand{\IntHeight}{8mm}
\renewcommand{\HorDist}{4.8mm}
\renewcommand{\VerDist}{4.8mm}
\renewcommand{\DblVerDist}{2*\VerDist} 
\renewcommand{\DblHorDist}{2*\HorDist}
\tikzstyle{oneaxis}=[draw,rectangle,thick,rounded corners = 1mm,minimum height=\IntHeight,adjust width=\IntAspect,inner sep=0,fill=white]
\fqwires[7*\HorDist]{w3}{4}
\anode[\VerDist/2]{node3}{w341}{}
\roaslantedlayer[black][\VerDist and \HorDist][\HorDist]{l3}{node3}{}{3}
\roagate[black][\DblHorDist]{g3}{l32}{}
\roagate[gray,dotted][\DblHorDist]{gb}{g3}{}
\coordinate[right= \HorDist/2 of g3] (a3);
\coordinate[right= \HorDist of a3] (a4);
\draw [myarrow] (a3) -- (a4);
\end{tikzpicture}
}\hfill%
\subfloat[]{%
\label{fig:gtl4}
\begin{tikzpicture}[on grid,auto]
\renewcommand{\IntHeight}{8mm}
\renewcommand{\HorDist}{4.8mm}
\renewcommand{\VerDist}{4.8mm}
\renewcommand{\DblVerDist}{2*\VerDist} 
\renewcommand{\DblHorDist}{2*\HorDist}
\tikzstyle{oneaxis}=[draw,rectangle,thick,rounded corners = 1mm,minimum height=\IntHeight,adjust width=\IntAspect,inner sep=0,fill=white]
\fqwires[5.5*\HorDist]{w4}{4}
\anode[\VerDist/2]{node4}{w441}{}
\roaslantedlayer[black][\VerDist and \HorDist][\HorDist]{l4}{node4}{}{3}
\roagate[black][2.5*\HorDist]{g4}{l42}{}
\end{tikzpicture}
}\\

\subfloat[]{%
\label{fig:gtl5}
\begin{tikzpicture}[on grid,auto]
\renewcommand{\IntHeight}{8mm}
\renewcommand{\HorDist}{4.8mm}
\renewcommand{\VerDist}{4.8mm}
\renewcommand{\DblVerDist}{2*\VerDist} 
\renewcommand{\DblHorDist}{2*\HorDist}
\tikzstyle{oneaxis}=[draw,rectangle,thick,rounded corners = 1mm,minimum height=\IntHeight,adjust width=\IntAspect,inner sep=0,fill=white]
\fqwires[8*\HorDist]{w1}{4}
\anode[\VerDist/2]{node1}{w141}{}
\roaslantedlayer[black][\VerDist and \HorDist][3*\HorDist]{l1}{node1}{}{3}
\broagate[black][\VerDist/2 and \HorDist]{g1}{w111}{}
\roagate[gray,dotted][3*\HorDist]{g4}{l12}{}
\coordinate[right= \HorDist/2 of g1] (a1);
\coordinate[right= 4.5*\HorDist of w121] (a2);
\coordinate[left= \HorDist/2 of g4] (a3);
\draw[rounded corners=2pt,myarrow](a1) -| (a2) |- (a3);        
\end{tikzpicture}
}\hfill%
\subfloat[]{%
\label{fig:gtl6}
\begin{tikzpicture}[on grid,auto]
\renewcommand{\IntHeight}{8mm}
\renewcommand{\HorDist}{4.8mm}
\renewcommand{\VerDist}{4.8mm}
\renewcommand{\DblVerDist}{2*\VerDist} 
\renewcommand{\DblHorDist}{2*\HorDist}
\tikzstyle{oneaxis}=[draw,rectangle,thick,rounded corners = 1mm,minimum height=\IntHeight,adjust width=\IntAspect,inner sep=0,fill=white]
\fqwires[8*\HorDist]{w1}{4}
\anode[\VerDist/2]{node1}{w141}{}
\roaslantedlayer[black][\VerDist and \HorDist][3*\HorDist]{l1}{node1}{}{3}
\broagate[gray,dotted][\VerDist/2 and \HorDist]{g1}{w121}{}
\roagate[black][3*\HorDist]{g4}{l11}{}
\coordinate[left= \HorDist/2 of g4] (a1);
\coordinate[right= 3.5*\HorDist of w131] (a2);
\coordinate[right= \HorDist/2 of g1] (a3);
\draw[rounded corners=2pt,myarrow](a1) -| (a2) |- (a3);        
\end{tikzpicture}
}\hfill%
\subfloat[]{%
\label{fig:gtl7}
\begin{tikzpicture}[on grid,auto]
\renewcommand{\IntHeight}{8mm}
\renewcommand{\HorDist}{4.8mm}
\renewcommand{\VerDist}{4.8mm}
\renewcommand{\DblVerDist}{2*\VerDist} 
\renewcommand{\DblHorDist}{2*\HorDist}
\tikzstyle{oneaxis}=[draw,rectangle,thick,rounded corners = 1mm,minimum height=\IntHeight,adjust width=\IntAspect,inner sep=0,fill=white]
\fqwires[8*\HorDist]{w1}{4}
\broagate[black][\VerDist/2 and \HorDist]{g1}{w121}{}
\broagate{g2}{g1}{}
\broagate[black][\VerDist/2 and 3*\HorDist]{ga}{w111}{}
\broagate{gb}{ga}{}
\broagate{gc}{gb}{}
\roagate[gray,dotted][3*\HorDist]{g3}{gb}{}
\roagate[gray,dotted][3*\HorDist]{g4}{ga}{}
\coordinate[right= \HorDist/2 of g2] (a1);
\coordinate[right= 4.5*\HorDist of w131] (a2);
\coordinate[left= \HorDist/2 of g3] (a3);
\draw[rounded corners=2pt,myarrow](a1) -- node[above, near start,fill=white] {\textbf{1}} ++(5.5mm,0mm) -| (a2) |- (a3);        
\coordinate[right= \HorDist/2 of g1] (a4);
\coordinate[right= 3.5*\HorDist of w121] (a5);
\coordinate[left= \HorDist/2 of g4] (a6);
\draw[rounded corners=2pt,myarrow](a4) -- node[above, near start,fill=white] {\textbf{2}} ++(5.5mm,0mm) -| (a5) |- (a6);        
\end{tikzpicture}
}\\
\label{fig:gtl}
\caption{(a) A free gate is moved in from the left side of the ascending cascade of three gates, (b) this gives a $\VEE$-shaped pattern of gates that are ready for a turnover, (c) the turnover operation results in a free gate on the right side of the ascending cascade that has moved one position down, to achieve final configuration shown in (d). 
Steps (a)-(d) are summarized with the concise notation in (e). (f) A free gate can be moved from right
to left, and (g) the order of turnover operations matters when moving multiple gates through a cascade.}
\end{figure}

An analogous operation can be performed for a free gate and a \emph{descending cascade} of gates,
or for a free gate on the right side of an ascending (\Cref{fig:gtl6}) or descending
cascade of gates. If we move multiple gates through a descending or ascending cascade, 
the order of turnover operations matters as illustrated in \Cref{fig:gtl7}.

\subsection{Square and triangle circuits}
\label{sec:SqTr}
We define two types of circuits with fixed patterns of two-spin gates: a circuit with a \emph{square}
pattern of gates and with a \emph{triangle} pattern of gates.

\begin{definition}[Square Circuit]
\label{def:SqCirc}
A square circuit on $N$ spins has $N$ vertical layers that are alternatingly starting from the first and second
spin. 
\end{definition}

\begin{definition}[Triangle Circuit]
\label{def:TrCirc}
A triangle circuit on $N$ spins has $N-1$ descending cascades with the $i$th descending 
cascade starting from spin $N-i$ and containing $i$ two-spin gates.
\end{definition}

\begin{figure}[hbtp]
\centering
~\hfill%
\begin{minipage}{0.45\textwidth}\centering
\subfloat[Square for $5$ spins]{%
\label{fig:square5}
\begin{tikzpicture}[on grid, auto]
\renewcommand{\IntHeight}{8mm}
\renewcommand{\HorDist}{4.8mm}
\renewcommand{\VerDist}{4.8mm}
\renewcommand{\DblVerDist}{2*\VerDist} 
\renewcommand{\DblHorDist}{2*\HorDist}
\tikzstyle{oneaxis}=[draw,rectangle,thick,rounded corners = 1mm,minimum height=\IntHeight,adjust width=\IntAspect,inner sep=0,fill=white]
\fqwires[6*\HorDist]{w1}{5}
\broalayer[black][\DblVerDist][\VerDist/2 and \HorDist]{l1}{w11}{}{2}
\broalayer{l2}{l11}{}{2}
\aroalayer{l3}{l21}{}{2}
\broalayer{l4}{l31}{}{2}
\aroalayer{l5}{l41}{}{2}
\end{tikzpicture}
}\end{minipage}\hfill%
\begin{minipage}{0.45\textwidth}\centering
\subfloat[Triangle for $5$ spins]{%
\label{fig:triangle5}
\begin{tikzpicture}[on grid, auto]
\renewcommand{\IntHeight}{8mm}
\renewcommand{\HorDist}{4.8mm}
\renewcommand{\VerDist}{4.8mm}
\renewcommand{\DblVerDist}{2*\VerDist} 
\renewcommand{\DblHorDist}{2*\HorDist}
\tikzstyle{oneaxis}=[draw,rectangle,thick,rounded corners = 1mm,minimum height=\IntHeight,adjust width=\IntAspect,inner sep=0,fill=white]
\fqwires[8*\HorDist]{w2}{5}
\anode[0.5*\HorDist]{a1}{w251}{}
\roapyramid[black][\HorDist]{py1}{a1}{}{4}
\end{tikzpicture}
}\end{minipage}\hfill~\\[10pt]

~\hfill%
\begin{minipage}{0.45\textwidth}\centering
\subfloat[Square for $6$ spins]{%
\label{fig:square}
\begin{tikzpicture}[on grid,auto]
\renewcommand{\IntHeight}{8mm}
\renewcommand{\HorDist}{4.8mm}
\renewcommand{\VerDist}{4.8mm}
\renewcommand{\DblVerDist}{2*\VerDist} 
\renewcommand{\DblHorDist}{2*\HorDist}
\tikzstyle{oneaxis}=[draw,rectangle,thick,rounded corners = 1mm,minimum height=\IntHeight,adjust width=\IntAspect,inner sep=0,fill=white]
\fqwires[7*\HorDist]{w1}{6}
\broalayer[black][\DblVerDist][\VerDist/2 and \HorDist]{l1}{w11}{}{3}
\broalayer{l2}{l11}{}{2}
\aroalayer{l3}{l21}{}{3}
\broalayer{l4}{l31}{}{2}
\aroalayer{l5}{l41}{}{3}
\broalayer{l6}{l51}{}{2}
\end{tikzpicture}
}\end{minipage}\hfill%
\begin{minipage}{0.45\textwidth}\centering
\subfloat[Square\,$\to$\,Triangle, step 1]{%
\label{fig:squarespace}
\begin{tikzpicture}[on grid,auto]
\renewcommand{\IntHeight}{8mm}
\renewcommand{\HorDist}{4.8mm}
\renewcommand{\VerDist}{4.8mm}
\renewcommand{\DblVerDist}{2*\VerDist} 
\renewcommand{\DblHorDist}{2*\HorDist}
\tikzstyle{oneaxis}=[draw,rectangle,thick,rounded corners = 1mm,minimum height=\IntHeight,adjust width=\IntAspect,inner sep=0,fill=white]
\fqwires[10*\HorDist]{w2}{6}
\broagate[black][\VerDist/2 and \HorDist]{g21}{w211}{}
\boaslantedlayer{g22}{g21}{}{3}
\boaslantedlayer[black][\VerDist and \HorDist][3*\VerDist]{g23}{g222}{}{5}
\roaslantedlayer{g24}{g231}{}{4}
\roaslantedlayer[black][\VerDist and \HorDist][2*\DblHorDist]{g25}{g241}{}{2}
\end{tikzpicture}
}\end{minipage}\hfill~\\[10pt]

~\hfill%
\begin{minipage}{0.45\textwidth}\centering
\subfloat[Square\,$\to$\,Triangle, step 2]{%
\label{fig:squarearrow}
\begin{tikzpicture}[on grid,auto]
\renewcommand{\IntHeight}{8mm}
\renewcommand{\HorDist}{4.8mm}
\renewcommand{\VerDist}{4.8mm}
\renewcommand{\DblVerDist}{2*\VerDist} 
\renewcommand{\DblHorDist}{2*\HorDist}
\tikzstyle{oneaxis}=[draw,rectangle,thick,rounded corners = 1mm,minimum height=\IntHeight,adjust width=\IntAspect,inner sep=0,fill=white]
\fqwires[10*\HorDist]{w3}{6}
\broagate[black][\VerDist/2 and \HorDist]{g11}{w311}{}
\boaslantedlayer{g12}{g11}{}{3}
\boaslantedlayer[black][\VerDist and \HorDist][3*\VerDist]{g13}{g122}{}{5}
\roaslantedlayer{g14}{g131}{}{4}
\roaslantedlayer[gray,dotted]{g16}{g141}{}{3}
\roaslantedlayer[black][\VerDist and \HorDist][2*\DblHorDist]{g15}{g141}{}{2}
\coordinate[right= \HorDist/2 of g123] (c11);
\coordinate[right= 5.5*\HorDist of w321] (c12);
\coordinate[right= \HorDist of g134] (c13);
\coordinate[right= 6.5*\HorDist of w331] (c14);
\coordinate[left= \HorDist/2 of g163] (c15);
\draw[rounded corners=2pt,myarrow](c11) -- node[below, near start,fill=white] {\textbf{1}} ++(5.5mm,0mm) -| (c12) |- (c13) -| (c14) |- (c15); 
\coordinate[right= \HorDist/2 of g122] (c21);
\coordinate[right= 4.5*\HorDist of w331] (c22);
\coordinate[right= \HorDist of g133] (c23);
\coordinate[right= 5.5*\HorDist of w341] (c24);
\coordinate[left= \HorDist/2 of g162] (c25);
\draw[rounded corners=2pt,myarrow](c21) -- node[below, near start,fill=white] {\textbf{2}} ++(5.5mm,0mm) -| (c22) |- (c23) -| (c24) |- (c25); 
\coordinate[right= \HorDist/2 of g121] (c31);
\coordinate[right= 3.5*\HorDist of w341] (c32);
\coordinate[right= \HorDist of g132] (c33);
\coordinate[right= 4.5*\HorDist of w351] (c34);
\coordinate[left= \HorDist/2 of g161] (c35);
\draw[rounded corners=2pt,myarrow](c31) -- node[below, near start,fill=white] {\textbf{3}} ++(5.5mm,0mm) -| (c32) |- (c33) -| (c34) |- (c35); 
\end{tikzpicture}
}\end{minipage}\hfill%
\begin{minipage}{0.45\textwidth}\centering
\subfloat[Triangle for $6$ spins]{%
\label{fig:almosttriangle}
\begin{tikzpicture}[on grid,auto]
\renewcommand{\IntHeight}{8mm}
\renewcommand{\HorDist}{4.8mm}
\renewcommand{\VerDist}{4.8mm}
\renewcommand{\DblVerDist}{2*\VerDist} 
\renewcommand{\DblHorDist}{2*\HorDist}
\tikzstyle{oneaxis}=[draw,rectangle,thick,rounded corners = 1mm,minimum height=\IntHeight,adjust width=\IntAspect,inner sep=0,fill=white]
\fqwires[10*\HorDist]{w4}{6}
\broagate[black][\VerDist/2 and \HorDist]{g21}{w411}{}
\boaslantedlayer[black][\VerDist and \HorDist][2*\DblVerDist]{g22}{g21}{}{5}
\roaslantedlayer{g23}{g221}{}{4}
\roaslantedlayer{g24}{g231}{}{3}
\roaslantedlayer{g25}{g241}{}{2}
\roagate[gray,dotted][\DblVerDist]{g26}{g251}{}
\coordinate[right= \HorDist/2 of g21] (cb1);
\coordinate[right= 4.5*\HorDist of w421] (cb2);
\coordinate[right= \HorDist of g224] (cb3);
\coordinate[right= 5.5*\HorDist of w431] (cb4);
\coordinate[right= \HorDist of g233] (cb5);
\coordinate[right= 6.5*\HorDist of w441] (cb6);
\coordinate[right= \HorDist of g242] (cb7);
\coordinate[right= 7.5*\HorDist of w451] (cb8);
\coordinate[left= \HorDist/2 of g26] (cb9);
\draw[rounded corners=2pt,myarrow](cb1) -- node[below, near start,fill=white] {\textbf{4}} ++(5.5mm,0mm) -| (cb2) |- (cb3) -| (cb4) |- (cb5) -| (cb6) |- (cb7) -| (cb8) |- (cb9);
\end{tikzpicture}
}\end{minipage}\hfill~\\
\label{fig:sqtr}
\caption{(a) Square and (b) triangle circuits for systems with $5$ spins. (c) A square circuit on $6$ spins and (d) the same square circuit with additional open space. 
(e) The three gates in the second ascending cascade of the square circuit are moved over to the bottom
half using a total of $6$ turnover operations, (f) finally the first gate is moved from the top left to 
bottom right using $4$ turnover operations resulting in a triangle circuit for $6$ spins.}
\end{figure}
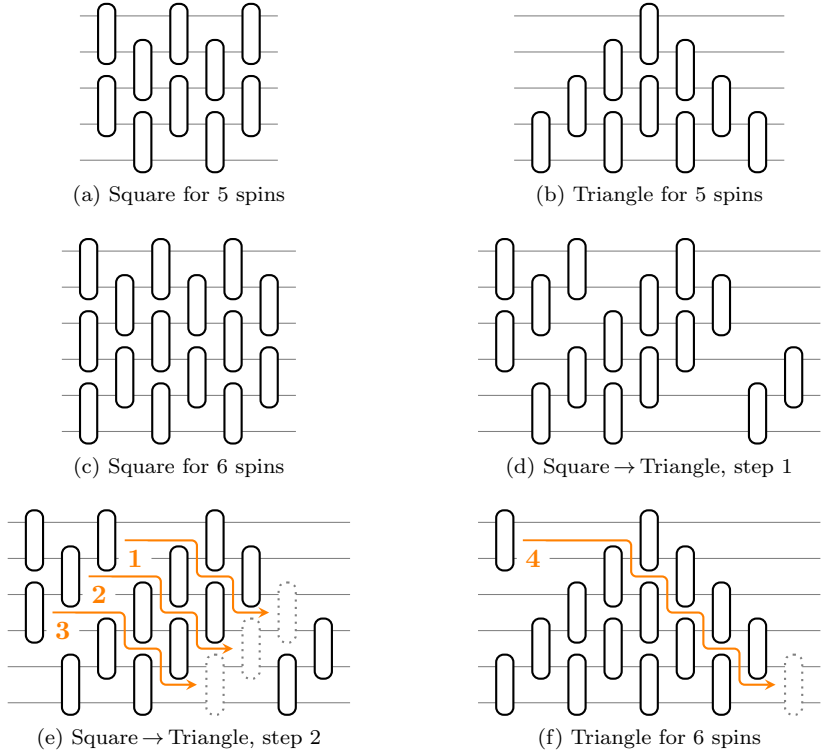

\Cref{fig:sqtr} illustrates square and triangular circuits for the odd and even
number of spins.
The number of gates in a square circuit is equal to the number of gates in a triangle pattern
and scales quadratically as $N(N-1)/2$.

Using the turnover operation, we can transform a square circuit
into a triangle circuit. The algorithm for $6$ spins is summarized in \Cref{fig:sqtr} 
and generalizes to any even number of spins.
In \Cref{fig:squarespace} we create sufficient open space in the square circuit such that we can use the turnover operation
to move the top half of gates in the square circuit over to the bottom half in \Cref{fig:squarearrow}, thereby creating
the triangle circuit in \Cref{fig:almosttriangle}.
With minor alterations the algorithm to convert a square to a triangle circuit can be extended
to systems with an odd number of spins.
Furthermore, all the turnover operations are reversible, which means that we can easily
transform a triangle circuit back to a square circuit by reversing the algorithm.
The computational complexity to go from a square to a triangle circuit on $N$ spins or vice versa 
scales as $\bigO(N^3)$.

\subsection{Merging gates with a triangle circuit}
\label{sec:Merge}

\begin{figure}[hbtp]
\centering
\begin{minipage}{0.52\textwidth}\centering
\subfloat[]{
\label{fig:slm1}
\begin{tikzpicture}[on grid,auto]
\renewcommand{\IntHeight}{8mm}
\renewcommand{\HorDist}{4.8mm}
\renewcommand{\VerDist}{4.8mm}
\renewcommand{\DblVerDist}{2*\VerDist} 
\renewcommand{\DblHorDist}{2*\HorDist}
\tikzstyle{oneaxis}=[draw,rectangle,thick,rounded corners = 1mm,minimum height=\IntHeight,adjust width=\IntAspect,inner sep=0,fill=white]
\fqwires[11*\HorDist]{w1}{5}
\anode[\VerDist/2]{node1}{w151}{}
\roaslantedlayer[black][\VerDist and \HorDist][\HorDist]{ll}{node1}{}{4}
\roapyramid[black][3*\HorDist]{py1}{ll1}{}{4}
\coordinate[right= \HorDist/2 of ll4] (ca1);
\coordinate[right= 6.5*\HorDist of w121] (ca2);
\coordinate[right= \HorDist of py113] (ca3);
\coordinate[right= 7.5*\HorDist of w131] (ca4);
\coordinate[right= \HorDist of py122] (ca5);
\coordinate[right= 8.5*\HorDist of w141] (ca6);
\coordinate[left= \HorDist/2 of py141] (ca7);
\draw[rounded corners=2pt,myarrow](ca1)  -- node[below, near start,fill=white] {\textbf{1}} ++(5.5mm,0mm) -| (ca2) |- (ca3) -| (ca4) |- (ca5) -| (ca6) |- (ca7);
\coordinate[right= \HorDist/2 of ll3] (cb1);
\coordinate[right= 5.5*\HorDist of w131] (cb2);
\coordinate[right= \HorDist of py112] (cb3);
\coordinate[right= 6.5*\HorDist of w141] (cb4);
\coordinate[left= \HorDist/2 of py131] (cb5);
\draw[rounded corners=2pt,myarrow](cb1)  -- node[below, near start,fill=white] {\textbf{2}} ++(5.5mm,0mm) -| (cb2) |- (cb3) -| (cb4) |- (cb5);
\coordinate[right= \HorDist/2 of ll2] (cc1);
\coordinate[right= 4.5*\HorDist of w141] (cc2);
\coordinate[left= \HorDist/2 of py121] (cc3);
\draw[rounded corners=2pt,myarrow](cc1)  -- node[below, near start,fill=white] {\textbf{3}} ++(5.5mm,0mm) -| (cc2) |- (cc3);
\coordinate[right= \HorDist/2 of ll1] (cd1);
\coordinate[left= \HorDist/2 of py111] (cd2);
\draw[rounded corners=2pt,myarrow](cd1)  -- node[below, near start,fill=white] {\textbf{4}} ++(5.5mm,0mm) -- (cd2);
\end{tikzpicture}
}\end{minipage}%
\begin{minipage}{0.48\textwidth}\centering
\subfloat[]{
\label{fig:slm2}
\begin{tikzpicture}[on grid,auto]
\renewcommand{\IntHeight}{8mm}
\renewcommand{\HorDist}{4.8mm}
\renewcommand{\VerDist}{4.8mm}
\renewcommand{\DblVerDist}{2*\VerDist} 
\renewcommand{\DblHorDist}{2*\HorDist}
\tikzstyle{oneaxis}=[draw,rectangle,thick,rounded corners = 1mm,minimum height=\IntHeight,adjust width=\IntAspect,inner sep=0,fill=white]
\fqwires[11*\HorDist]{w2}{5}
\anode[\VerDist/2]{node2}{w251}{}
\roaslantedlayer[black][\VerDist and \HorDist][\HorDist]{ll}{node2}{}{4}
\roapyramid[black][\HorDist]{py2}{node2}{}{4}
\roagate[black][3*\HorDist]{g1}{py214}{}
\broagate{g2}{g1}{}
\broagate{g3}{g2}{}
\broagate{g4}{g3}{}
\coordinate[left= \HorDist/2 of g1] (da1);
\coordinate[left= 6.5*\HorDist of w222] (da2);
\coordinate[left= \HorDist of py223] (da3);
\coordinate[left= 7.5*\HorDist of w232] (da4);
\coordinate[left= \HorDist of py222] (da5);
\coordinate[left= 8.5*\HorDist of w242] (da6);
\coordinate[right= \HorDist/2 of py211] (da7);
\draw[rounded corners=2pt,myarrow](da1) -- node[below, near start,fill=white] {\textbf{1}} ++(-5.5mm,0mm) -| (da2) |- (da3) -| (da4) |- (da5) -| (da6) |- (da7);
\coordinate[left= \HorDist/2 of g2] (db1);
\coordinate[left= 5.5*\HorDist of w232] (db2);
\coordinate[left= \HorDist of py232] (db3);
\coordinate[left= 6.5*\HorDist of w242] (db4);
\coordinate[right= \HorDist/2 of py221] (db5);
\draw[rounded corners=2pt,myarrow](db1)  -- node[below, near start,fill=white] {\textbf{2}} ++(-5.5mm,0mm) -| (db2) |- (db3) -| (db4) |- (db5);
\coordinate[left= \HorDist/2 of g3] (dc1);
\coordinate[left= 4.5*\HorDist of w242] (dc2);
\coordinate[right= \HorDist/2 of py231] (dc3);
\draw[rounded corners=2pt,myarrow](dc1)  -- node[below, near start,fill=white] {\textbf{3}} ++(-5.5mm,0mm) -| (dc2) |- (dc3);
\coordinate[left= \HorDist/2 of g4] (dd1);
\coordinate[right= \HorDist/2 of py241] (dd2);
\draw[rounded corners=2pt,myarrow](dd1)  -- node[below, near start,fill=white] {\textbf{4}} ++(-5.5mm,0mm) -- (dd2);
\end{tikzpicture}
}\end{minipage}\\[10pt]

\begin{minipage}{0.52\textwidth}\centering
\subfloat[]{
\label{fig:vert1}
\begin{tikzpicture}[on grid,auto]
\renewcommand{\IntHeight}{8mm}
\renewcommand{\HorDist}{4.8mm}
\renewcommand{\VerDist}{4.8mm}
\renewcommand{\DblVerDist}{2*\VerDist} 
\renewcommand{\DblHorDist}{2*\HorDist}
\tikzstyle{oneaxis}=[draw,rectangle,thick,rounded corners = 1mm,minimum height=\IntHeight,adjust width=\IntAspect,inner sep=0,fill=white]

\fqwires[13*\HorDist]{w1}{6}
\anode[\VerDist/2]{node1}{w161}{}
\roapyramid[black][\HorDist]{py1}{node1}{}{5}
\roalayer[black][\DblVerDist][6*\HorDist]{l1}{py115}{}{3}
\roalayer[black][\DblVerDist][6*\HorDist]{l2}{py124}{}{2}
\coordinate[left= \HorDist/2 of l11] (da1);
\coordinate[left= 7.5*\HorDist of w122] (da2);
\coordinate[left= \HorDist of py124] (da3);
\coordinate[left= 8.5*\HorDist of w132] (da4);
\coordinate[left= \HorDist of py123] (da5);
\coordinate[left= 9.5*\HorDist of w142] (da6);
\coordinate[left= \HorDist of py122] (da7);
\coordinate[left= 10.5*\HorDist of w152] (da8);
\coordinate[right= \HorDist/2 of py111] (da9);
\draw[rounded corners=2pt,myarrow](da1) -- node[below, near start,fill=white] {\textbf{1}} ++(-5.5mm,0mm) -| (da2) |- (da3) -| (da4) |- (da5) -| (da6) |- (da7) -| (da8) |- (da9);
\coordinate[left= \HorDist/2 of l12] (db1);
\coordinate[left= 5.5*\HorDist of w142] (db2);
\coordinate[left= \HorDist of py142] (db3);
\coordinate[left= 6.5*\HorDist of w152] (db4);
\coordinate[right= \HorDist/2 of py131] (db5);
\draw[rounded corners=2pt,myarrow](db1) -- node[below, near start,fill=white] {\textbf{1}} ++(-5.5mm,0mm) -| (db2) |- (db3) -| (db4) |- (db5);
\coordinate[left= \HorDist/2 of l13] (dc1);
\coordinate[right= \HorDist/2 of py151] (dc2);
\draw[rounded corners=2pt,myarrow](dc1) -- node[below, near start,fill=white] {\textbf{1}} ++(-5.5mm,0mm) -- (dc2);
\bnode[1mm]{mn1}{w161}{}
\end{tikzpicture}
}\end{minipage}%
\begin{minipage}{0.48\textwidth}\centering
\subfloat[]{
\label{fig:vert2}
\begin{tikzpicture}[on grid,auto]
\renewcommand{\IntHeight}{8mm}
\renewcommand{\HorDist}{4.8mm}
\renewcommand{\VerDist}{4.8mm}
\renewcommand{\DblVerDist}{2*\VerDist} 
\renewcommand{\DblHorDist}{2*\HorDist}
\tikzstyle{oneaxis}=[draw,rectangle,thick,rounded corners = 1mm,minimum height=\IntHeight,adjust width=\IntAspect,inner sep=0,fill=white]
\fqwires[12*\HorDist]{w2}{6}
\anode[\VerDist/2]{node2}{w261}{}
\roapyramid[black][\HorDist]{py2}{node2}{}{5}
\roalayer[black][\DblVerDist][5*\HorDist]{l3}{py224}{}{2}
\coordinate[left= \HorDist/2 of l31] (dd1);
\coordinate[left= 5.5*\HorDist of w232] (dd2);
\coordinate[left= \HorDist of py233] (dd3);
\coordinate[left= 6.5*\HorDist of w242] (dd4);
\coordinate[left= \HorDist of py232] (dd5);
\coordinate[left= 7.5*\HorDist of w252] (dd6);
\coordinate[right= \HorDist/2 of py221] (dd7);
\draw[rounded corners=2pt,myarrow](dd1) -- node[below, near start,fill=white] {\textbf{2}} ++(-5.5mm,0mm) -| (dd2) |- (dd3) -| (dd4) |- (dd5) -| (dd6) |- (dd7);
\coordinate[left= \HorDist/2 of l32] (de1);
\coordinate[left= 3.5*\HorDist of w252] (de2);
\coordinate[right= \HorDist/2 of py241] (de3);
\draw[rounded corners=2pt,myarrow](de1) -- node[below, near start,fill=white] {\textbf{2}} ++(-5.5mm,0mm) -| (de2) |- (de3);
\bnode[1mm]{mn2}{w261}{}
\end{tikzpicture}
}\end{minipage}\\
\label{fig:merge}
\caption{(a) Sequential annihilation of an ascending sequence of gates positioned 
to the left of a $5$-spin triangle circuit through repeated turnover and fuse operations, (b)
similar for a sequence of gates on the right of a $5$-spin triangle circuit.
A vertical layer of gates can be merged in parallel with a $6$-spin triangle circuit by first
merging all gates acting on odd spins (c) and afterwards merging the remaining gates on the even spins (d).}
\end{figure}
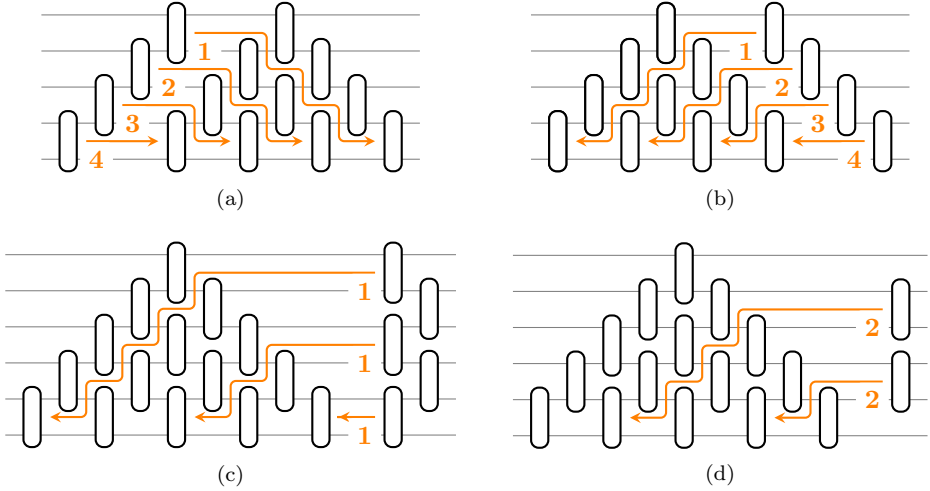

We will now proceed to show that a triangle is the minimal circuit for two-spin
gates that allow for a turnover and fusion operation.
This holds because every two-spin gate on either the left or right side of a triangle circuit
can be repeatedly turned over with gates in the triangle circuit until it eventually can be
fused with a gate at the bottom of the triangle circuit.
This is illustrated in \Cref{fig:slm1,fig:slm2} for a $5$-spin triangle circuit with a cascade of 
gates on both the left and right side.

For a Trotter-based implementation of a simulation circuit, the gates often come in vertical
layers that act alternating on even and odd spins. In that case the annihilation of the gates
acting on even spins can all be done in parallel, and the same is true for the gates acting on
odd spins. This is shown in \Cref{fig:vert1,fig:vert2}.

The average computational cost of merging a single gate with a triangle circuit scales as $\bigO(N)$,
and the cost of merging a complete vertical layer of gates with a triangle is $\bigO(N^2)$.
There is no difference in computational complexity between disordered and ordered Hamiltonians,
both require $\bigO(n_t N^2)$ operations to annihilate $n_t$ time-step.
For time-independent Hamiltonians, we can again reduce the complexity to $\bigO(\log_2(n_t) N^2)$
by using the merging algorithm recursively.

\subsection{Circuit compression algorithms}

The complete circuit compression algorithms are outlined in \Cref{alg:compr1}
for time-dependent Hamiltonians and in \Cref{alg:compr2} for the time-independent
case.
The input to \Cref{alg:compr1} is a Trotter circuit $C$ on $N$ spins
with $n_t$ time-steps or $2n_t$ vertical layers. We assume that $n_t > N/2$
as otherwise minimal depth is not reached.
In line~\ref{line:toTr}, we take the first $N$ vertical layers out of $C$ and transform
it to triangle representation $C^{\prime}$.
The for loop runs over the remaining layers and merges them into $C^{\prime}$
using the approach of \Cref{fig:vert1,fig:vert2}.
In the end $C^{\prime}$ is transformed back to a square circuit of depth $N$
that is equivalent to the input circuit $C$ with depth $2 n_t$.

\begin{algorithm}[hbtp]
\SetKwFunction{TriangleCircuit}{TriangleCircuit}
\SetKwFunction{SquareCircuit}{SquareCircuit}
\SetKwFunction{Merge}{Merge}
\SetKwFunction{MergeLayer}{MergeLayer}
\KwData{Trotter circuit $C$ on $N$ spins with $n_t$ time-steps, $n_t > N/2$.}
\KwResult{Compressed $N \times N$ Trotter circuit $C^{\prime}$ equivalent to $C$.}
\BlankLine
$C^{\prime} \longleftarrow$ \TriangleCircuit{$C[:,1{:}N]$}\label{line:toTr}

\For{$l\leftarrow N+1$ \KwTo $2 n_t$}{
	\MergeLayer{$C^{\prime}, C[:,l]$}
}
$C^{\prime} \longleftarrow$ \SquareCircuit{$C^{\prime}$}
\caption{Compression of time-dependent Hamiltonians}
\label{alg:compr1}
\end{algorithm}

In the time-independent case, we first create a minimal depth square $C^{\prime}$ from
the fixed time-step $C_{TS}$ in line~\ref{line:mds}.
This square is again transformed to a triangle that is repeatedly merged with itself 
in line~\ref{line:rmerge} by splitting the triangle into cascades and using the method from
\Cref{fig:slm2}. This doubles the total simulation time in every loop.

Both algorithms end with a conversion to square format as this form has a shallower 
circuit depth compared to the triangle format. This makes square circuits better suited for NISQ devices,
even though both circuits have the same number of gates.

\begin{algorithm}[hbtp]
\setcounter{AlgoLine}{0}
\SetKwFunction{TriangleCircuit}{TriangleCircuit}
\SetKwFunction{SquareCircuit}{SquareCircuit}
\SetKwFunction{Merge}{Merge}
\SetKwFunction{MergeTriangle}{MergeTriangle}
\KwData{Trotter circuit $C_{TS}$ for a single time-step, total number of time-steps $n_t = 2^{\tau} N/2$.}
\KwResult{$N \times N$ Trotter circuit $C^{\prime}$ for $n_t$ time-steps.}
\BlankLine
\tcp*[h]{Fill $C^{\prime}$ to square form with $N/2$ time-steps $C_{TS}$.}

\For{$l \leftarrow 1$ \KwTo $N/2$}{
	$C^{\prime}[:,l{:}l{+}1] \longleftarrow C_{TS}$ \label{line:mds}
}
$C^{\prime} \longleftarrow$ \TriangleCircuit{$C^{\prime}$}

\tcp*[h]{Repeatedly merge $C^{\prime}$ with itself up to total time-steps of $2^{\tau} N/2$.}

\For{$t \leftarrow 1$ \KwTo $\tau$}{
	\MergeTriangle{$C^{\prime},C^{\prime}$}\label{line:rmerge}
}
$C^{\prime} \longleftarrow$ \SquareCircuit{$C^{\prime}$}
\caption{Compression of time-independent Hamiltonians}
\label{alg:compr2}
\end{algorithm}


\section{Kitaev chains and XY models}
\label{sec:KXY}
In this section we show that two classes of closely related Hamiltonians, 
respectively known as Kitaev chains and XY models, can be implemented
with gates that satisfy the fusion and turnover conditions introduced in the
previous section.
It follows that the discretized time-evolution operator \eqref{eq:dtimeevol}
can be implemented in a quantum circuit with depth $\bigO(N)$ for these models.

\subsection{Kitaev chains}
\label{sec:Kitaev}

\begin{figure}[hbtp]
\centering
~\hfill%
\subfloat[]{%
\label{fig:Kitaev1}
\begin{tikzpicture}[on grid, auto]
\fqwires[3*\HorDist]{w1}{5}
\broagate[myred][\VerDist/2 and \HorDist]{g1}{w111}{$y$}
\boagate[mygreen][\DblVerDist]{g2}{g1}{$z$}
\broagate[myblue][\VerDist and \HorDist]{g3}{g1}{$x$}
\broagate[myblue][\VerDist and \HorDist]{g4}{g2}{$x$}
\end{tikzpicture}
}\hfill%
\subfloat[]{%
\label{fig:Kitaev2}
\begin{tikzpicture}[on grid, auto]
\fqwires[6*\HorDist]{w2}{5}
\broagate[myred][\VerDist/2 and \HorDist]{g5}{w211}{$y$}
\boagate[mygreen][\DblVerDist]{g6}{g5}{$z$}
\broagate[myblue][\VerDist and \HorDist]{g7}{g5}{$x$}
\broagate[myblue][\VerDist and \HorDist]{g8}{g6}{$x$}
\aroagate[myred][\VerDist and \HorDist]{g9}{g7}{$y$}
\aroagate[mygreen][\VerDist and \HorDist]{g10}{g8}{$z$}
\broagate[myblue][\VerDist and \HorDist]{g11}{g9}{$x$}
\broagate[myblue][\VerDist and \HorDist]{g12}{g10}{$x$}
\aroagate[myred][\VerDist and \HorDist]{g13}{g11}{$y$}
\aroagate[mygreen][\VerDist and \HorDist]{g14}{g12}{$z$}
\end{tikzpicture}
}\hfill%
\subfloat[]{%
\label{fig:Kitaev3}
\begin{tikzpicture}[on grid, auto]
\rqwires[8*\HorDist]{w3}{w212}{5}
\aroagate[myblue][\VerDist/2 and \HorDist]{g21}{w351}{$x$}
\aroagate[mygreen][\VerDist and \HorDist]{g22}{g21}{$z$}
\aroagate[myblue][\VerDist and \HorDist]{g23}{g22}{$x$}
\aroagate[myred][\VerDist and \HorDist]{g24}{g23}{$y$}
\roagate[myblue][\DblHorDist]{g25}{g21}{$x$}
\aroagate[mygreen][\VerDist and \HorDist]{g26}{g25}{$z$}
\aroagate[myblue][\VerDist and \HorDist]{g27}{g26}{$x$}
\roagate[myblue][\DblHorDist]{g28}{g25}{$x$}
\aroagate[mygreen][\VerDist and \HorDist]{g29}{g28}{$z$}
\roagate[myblue][\DblHorDist]{g30}{g28}{$x$}
\end{tikzpicture}
}\hfill~\\
\label{fig:Kitaev}
\caption{(a) A single time-step for the Kitaev chain \eqref{eq:Kitaevex}, (b)
a square circuit for this Kitaev chain, and (c) a triangle circuit for
this Kitaev chain.}
\end{figure}
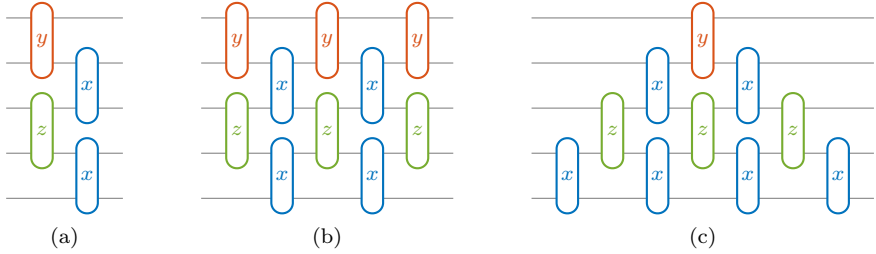

A Kitaev chain is a Hamiltonian of the form
\begin{align}
\label{eq:HKitaev}
H(t) & = \sum_{i=1}^{N-1} J^{\alpha_{i}}_i(t) \, \sigma^{\alpha_{i}}_i \sigma^{\alpha_{i}}_{i+1},
\end{align}
with the restriction that two neighboring spins $i$ and $i+1$ have different type of interaction, 
i.e.~$\alpha_i \neq \alpha_{i+1}$.
For example,
\begin{align}
\label{eq:Kitaevex}
H(t) = J^{y}_1(t) \, \sigma^{y}_1 \sigma^{y}_{2} \ + \ J^{x}_2(t) \, \sigma^{x}_2 \sigma^{x}_{3} 
\ + \ J^{z}_3(t) \, \sigma^{z}_3 \sigma^{z}_{4} \ + \ J^{x}_4(t) \, \sigma^{x}_4 \sigma^{x}_{5},
\end{align}
is a $5$-spin Kitaev chain.
Using a Trotter decomposition in even and odd terms, we get
the following approximation for a single time-step:
\begin{align*}
U_{\tau}(\Delta t) & = \exp(-\I H(t) \Delta(t) ), \\
& = \exp(-\I H_{\text{even}}(t) \Delta(t) ) \, \exp(-\I H_{\text{odd}}(t) \Delta(t) ) + \bigO(\Delta t^2), \\
& \approx 
R^{xx}_{2}(2 J^{x}_2(t_{\tau}) \Delta t) \,
R^{xx}_{4}(2 J^{x}_4(t_{\tau}) \Delta t) \,
R^{yy}_{1}(2 J^{y}_1(t_{\tau}) \Delta t) \,
R^{zz}_{3}(2 J^{z}_3(t_{\tau}) \Delta t).
\end{align*}
This is illustrated in a circuit diagram on the
left side of \Cref{fig:Kitaev}.
A square circuit for the Hamiltonian \eqref{eq:Kitaevex} that satisfies \Cref{def:SqCirc} is 
shown in the middle of \Cref{fig:Kitaev}.
It consists of two complete time-steps and a single vertical layer of gates
acting on the odd spins.
Using the turnover operation from \Cref{lem:mato},
we can use the algorithm from \Cref{sec:SqTr} to transform this square circuit to the
equivalent triangle circuit shown on the right of \Cref{fig:Kitaev}.
Remark that the number of gates of each type is changed due to \Cref{lem:mato}, but the type of gate acting on every
pair of spins is preserved.
Afterwards, we can use the turnover and fusion operation from \Cref{lem:pfus}\ref{itm:fusaa}
to compress the circuit to constant depth by means of \Cref{alg:compr1} or \Cref{alg:compr2}.
This shows that we can always compress the circuit for the simulation of a Kitaev to square or triangle form.

\subsection{XY models}
\label{sec:XY}

\begin{figure}[hbtp]
\centering
~\hfill%
\begin{minipage}{0.45\textwidth}\centering
\subfloat[]{%
\label{fig:XY1}
\begin{tikzpicture}[on grid, auto]
\fqwires[8*\HorDist]{w1}{5}
\aroagate[myred][\VerDist/2 and \HorDist]{g11}{w151}{$y$}
\aroagate[myblue][\VerDist and \HorDist]{g12}{g11}{$x$}
\aroagate[myred][\VerDist and \HorDist]{g13}{g12}{$y$}
\aroagate[myblue][\VerDist and \HorDist]{g14}{g13}{$x$}
\roagate[myred][\DblHorDist]{g15}{g11}{$y$}
\aroagate[myblue][\VerDist and \HorDist]{g16}{g15}{$x$}
\aroagate[myred][\VerDist and \HorDist]{g17}{g16}{$y$}
\roagate[myred][\DblHorDist]{g18}{g15}{$y$}
\aroagate[myblue][\VerDist and \HorDist]{g19}{g18}{$x$}
\roagate[myred][\DblHorDist]{g20}{g18}{$y$}
\end{tikzpicture}
}\end{minipage}\hfill%
\begin{minipage}{0.45\textwidth}\centering
\subfloat[]{%
\label{fig:XY2}
\begin{tikzpicture}[on grid, auto]
\fqwires[8*\HorDist]{w2}{5}
\aroagate[myblue][\VerDist/2 and \HorDist]{g21}{w251}{$x$}
\aroagate[myred][\VerDist and \HorDist]{g22}{g21}{$y$}
\aroagate[myblue][\VerDist and \HorDist]{g23}{g22}{$x$}
\aroagate[myred][\VerDist and \HorDist]{g24}{g23}{$y$}
\roagate[myblue][\DblHorDist]{g25}{g21}{$x$}
\aroagate[myred][\VerDist and \HorDist]{g26}{g25}{$y$}
\aroagate[myblue][\VerDist and \HorDist]{g27}{g26}{$x$}
\roagate[myblue][\DblHorDist]{g28}{g25}{$x$}
\aroagate[myred][\VerDist and \HorDist]{g29}{g28}{$y$}
\roagate[myblue][\DblHorDist]{g30}{g28}{$x$}
\end{tikzpicture}
}\end{minipage}\hfill~\\[10pt]

~\hfill%
\begin{minipage}{0.45\textwidth}\centering
\subfloat[]{%
\label{fig:XY3}
\begin{tikzpicture}[on grid, auto]
\fqwires[8*\HorDist]{w3}{5}
\artagate[\VerDist/2 and \HorDist]{myblue,myred}{g21}{w351}{$x$}{$y$}
\artagate[\VerDist and \HorDist]{myblue,myred}{g22}{g21}{$x$}{$y$}
\artagate[\VerDist and \HorDist]{myblue,myred}{g23}{g22}{$x$}{$y$}
\artagate[\VerDist and \HorDist]{myblue,myred}{g24}{g23}{$x$}{$y$}
\rtagate[\DblHorDist]{myblue,myred}{g25}{g21}{$x$}{$y$}
\artagate[\VerDist and \HorDist]{myblue,myred}{g26}{g25}{$x$}{$y$}
\artagate[\VerDist and \HorDist]{myblue,myred}{g27}{g26}{$x$}{$y$}
\rtagate[\DblHorDist]{myblue,myred}{g28}{g25}{$x$}{$y$}
\artagate[\VerDist and \HorDist]{myblue,myred}{g29}{g28}{$x$}{$y$}
\rtagate[\DblHorDist]{myblue,myred}{g30}{g28}{$x$}{$y$}
\end{tikzpicture}
}\end{minipage}\hfill%
\begin{minipage}{0.45\textwidth}\centering
\subfloat[]{%
\label{fig:XY4}
\begin{tikzpicture}[on grid, auto]
\fqwires[2*\HorDist]{w2}{2}
\brtagate[\VerDist/2 and \HorDist]{myblue,mygreen}{g2}{w211}{$x$}{$z$}
\bqwires[2*\HorDist][\VerDist][2*\VerDist]{w3}{w221}{2}
\brtagate[\VerDist/2 and \HorDist]{myred,mygreen}{g3}{w311}{$y$}{$z$}
\end{tikzpicture}
}\end{minipage}\hfill~\\
\label{fig:XYTriangle}
\caption{Triangle circuits for the Kitaev chains $H_{XY}$ (a) and $H_{YX}$ (b).
(c) Combination of two Kitaev chains in a single triangle circuit for the complete XY model,
(d) two-spin gates for XZ and YZ models.}
\end{figure}
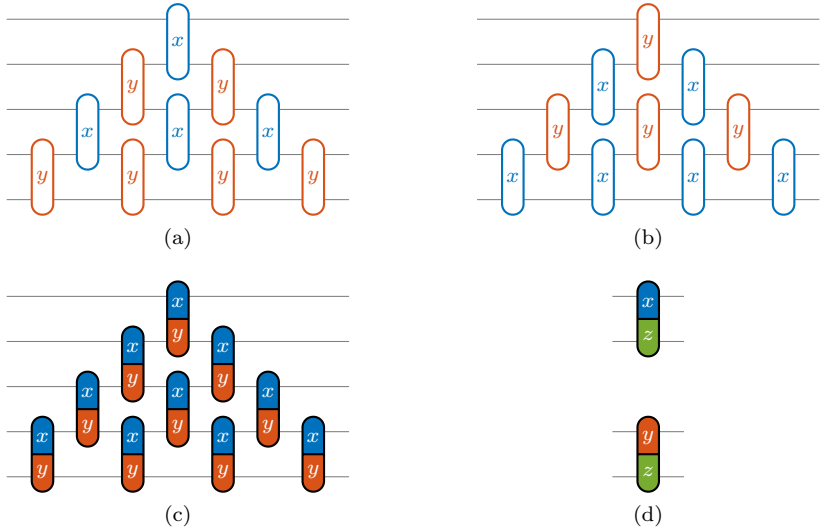

The class of XY Hamiltonians is given by
\begin{align}
\label{eq:XY}
H(t) & = \sum_{i=1}^{N-1} J^{\alpha}_i(t)  \, \sigma_{i}^{\alpha} \sigma_{i+1}^{\alpha} + 
J^{\beta}_i(t)  \, \sigma_{i}^{\beta} \sigma_{i+1}^{\beta},
&  \alpha, \beta \in \lbrace x, y, z \rbrace,
\end{align}
where $\alpha \neq \beta$ and all parameters $J^{\alpha}_i, J^{\beta}_i$ are nonzero.
This class consists of XY, XZ, and YZ Hamiltonians for the appropriate choices of
$\alpha$ and $\beta$.
If we split \eqref{eq:XY} in $H(t) = H_{\alpha\beta}(t) + H_{\beta\alpha}(t)$ with
\begin{align*}
H_{\alpha\beta}(t) & = \sum_{\text{odd } i} J^{\alpha}_i(t)  \, \sigma_{i}^{\alpha} \sigma_{i+1}^{\alpha} 
+ \sum_{\text{even } i} J^{\beta}_i(t)  \, \sigma_{i}^{\beta} \sigma_{i+1}^{\beta}, \\
H_{\beta\alpha}(t) & = \sum_{\text{odd } i} J^{\beta}_i(t)  \, \sigma_{i}^{\beta} \sigma_{i+1}^{\beta}
+ \sum_{\text{even } i} J^{\alpha}_i(t)  \, \sigma_{i}^{\alpha} \sigma_{i+1}^{\alpha},
\end{align*}
we see that we have rewritten it as a sum of two regular Kitaev chains.
Using the Trotter decomposition and circuit compression described in \Cref{sec:Kitaev}
for each of the two Kitaev chains $H_{\alpha\beta}(t)$ and $H_{\beta\alpha}(t)$,
we find that we can get a circuit for $H(t)$ as a product of two triangle circuits
for the Kitaev chain.
This is illustrated in the first row of \Cref{fig:XYTriangle} for the case of the XY model.
It follows from the commutativity properties \ref{itm:pcomaa} and \ref{itm:pcomaabb} 
of \Cref{lem:pcom} that we can combine these two Kitaev chains into one triangle circuit,
as shown on the second row of \Cref{fig:XYTriangle}.
The gates in this circuit are two-axes rotation gates defined in \Cref{lem:2axes},
and are a product of $R^{xx}$ and $R^{yy}$ rotations.
Because of the commutativity of the two Kitaev chains, they can be simulated separately.
However, for certain types of quantum hardware, a product of $R^{xx}$ and $R^{yy}$ rotations
can be evaluated at approximately the same cost as a single two-spin Pauli rotation~\cite{Bassman2021}.
In that case it makes sense to simulate the two Kitaev chains simultaneously.

The elementary gates for XY, XZ and YZ Hamiltonians are listed in \Cref{fig:XY3,fig:XY4}. Within each class,
these gates can be fused and turned over by separating them into their two-spin Pauli rotations
and using the results from \Cref{sec:props}.


\section{Transverse-field XY and Ising models}
\label{sec:transverse}

In this section, we discuss the most general class of Hamiltonians for which our compression
algorithm works. These are the TFXY Hamiltonians already introduced in \eqref{eq:TFXY}.
In \Cref{sec:TFXY} we briefly discuss the case of the full TFXY model, but more details about
the implementation of the fusion and turnover operation are deferred to \Cref{sec:implementation}.
\Cref{sec:TFIM} introduces the transverse-field Ising model (TFIM) as a special case of TFXY that
can be treated separately under some conditions.

\subsection{TFXY model}
\label{sec:TFXY}

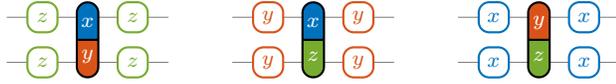
\begin{figure}[hbtp]
\centering
\begin{tikzpicture}[on grid, auto]
\fqwires[4*\HorDist]{w1}{2}
\rextgate[mygreen]{gate2}{w111}{$z$}
\rextgate[mygreen]{gate3}{w121}{$z$}
\brtagate[\VerDist/2 and 2*\HorDist]{myblue,myred}{gate1}{w111}{$x$}{$y$}
\lextgate[mygreen]{gate4}{w112}{$z$}
\lextgate[mygreen]{gate5}{w122}{$z$}
\rqwires[4*\HorDist]{w2}{w112}{2}
\rextgate[myred]{gate2}{w211}{$y$}
\rextgate[myred]{gate3}{w221}{$y$}
\brtagate[\VerDist/2 and 2*\HorDist]{myblue,mygreen}{gate1}{w211}{$x$}{$z$}
\lextgate[myred]{gate4}{w212}{$y$}
\lextgate[myred]{gate5}{w222}{$y$}
\rqwires[4*\HorDist]{w3}{w212}{2}
\rextgate[myblue]{gate2}{w311}{$x$}
\rextgate[myblue]{gate3}{w321}{$x$}
\brtagate[\VerDist/2 and 2*\HorDist]{myred,mygreen}{gate1}{w311}{$y$}{$z$}
\lextgate[myblue]{gate4}{w312}{$x$}
\lextgate[myblue]{gate5}{w322}{$x$}
\end{tikzpicture}
\caption{Two-spin gates for TFXY, TFXZ, and TFYZ Hamiltonians.}
\label{fig:TFXY}
\end{figure}

The two-spin gates for TFXY, TFXZ, and TFYZ Hamiltonians are shown in \Cref{fig:TFXY}.
The existence of the TFXY fusion and turnover operations is proven in \cite{TrotterCompression}
by considering the 
Hamiltonian algebras of 2- and 3-spin TFXY models and their Cartan decompositions.
It is further shown that the TFXY fusion requires two turnovers of Euler decompositions of $\SU2$,
\Cref{lem:turnoverSU2}, and that
a TFXY turnover can be done by invoking this lemma 32 times.
In our compression algorithms we use more efficient implementations of the fusion
and turnover operations described in \Cref{sec:implementation}.
With these two operations, we can use \Cref{alg:compr1,alg:compr2} to compress TFXY circuits.

\subsection{TFIM model}
\label{sec:TFIM}

The transverse-field Ising model is a special case of the TFXY Hamiltonian with
only one non-zero coupling term,
\begin{align}
H(t) = 
& \sum_{i=1}^{N-1} J^{\alpha}_i(t) \, \sigma_{i}^{\alpha} \sigma_{i+1}^{\alpha} 
\ + \ \sum_{i=1}^{N} h^{\beta}_i(t) \, \sigma_{i}^\beta,
&  \alpha, \beta \in \lbrace x, y, z \rbrace,
\label{eq:TFIM}
\end{align}
where $\alpha \neq \beta$.
The most frequently studied cases in the literature are $\alpha = x$, $\beta = z$, and
$\alpha = z$, $\beta = x$. We present our discussion for the former case, but all results remain valid for other
choices.

One straightforward approach for compressing a TFIM circuit is to use the more general TFXY gates
from \Cref{sec:TFXY} and set the $J^y$ parameters to zero.
While this leads to a valid compression algorithm, the $J^y$ parameters in the triangle circuit
will become nonzero throughout the procedure as the gate
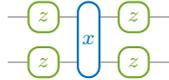
\begin{equation*}
\begin{tikzpicture}[on grid, auto]
\fqwires[4*\HorDist]{w1}{2}
\rextgate[mygreen]{gate2}{w111}{$z$}
\rextgate[mygreen]{gate3}{w121}{$z$}
\broagate[myblue][\VerDist/2 and 2*\HorDist]{gate1}{w111}{$x$}
\lextgate[mygreen]{gate4}{w112}{$z$}
\lextgate[mygreen]{gate5}{w122}{$z$}
\end{tikzpicture}
\vspace{-5pt}
\end{equation*}
doesn't allow for fusion and turnover operations.

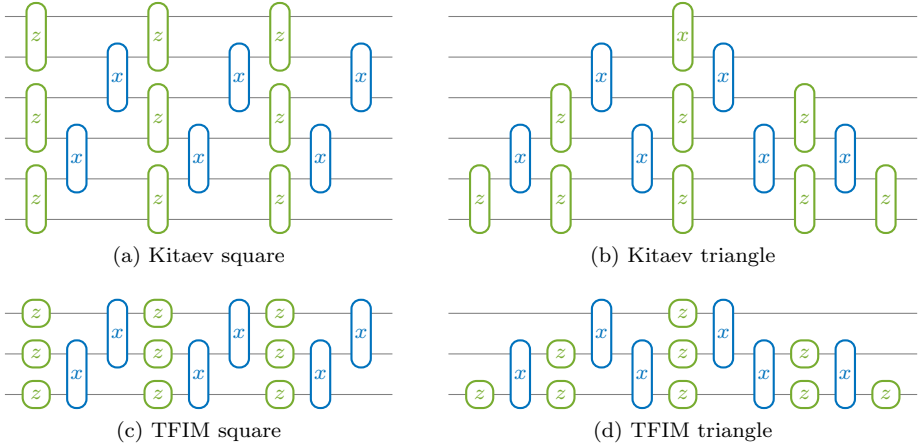
\begin{figure}[hbtp]
\centering
~\hfill%
\subfloat[Kitaev square]{%
\label{fig:TFK1}
\begin{tikzpicture}[on grid, auto]
\renewcommand{\ExtSize}{3.6mm} 
\renewcommand{\IntHeight}{9mm}
\renewcommand{\HorDist}{5.4mm}
\renewcommand{\VerDist}{5.4mm}
\renewcommand{\DblVerDist}{2*\VerDist} 
\renewcommand{\DblHorDist}{2*\HorDist}
\tikzstyle{oneaxis}=[draw,rectangle,thick,rounded corners = 1.2mm,minimum height=\IntHeight,adjust width=\IntAspect,inner sep=0,fill=white]
\fqwires[10*\HorDist]{w1}{6}
\broagate[mygreen][\VerDist/2 and \HorDist]{g1}{w111}{$z$}
\boagate[mygreen][\DblVerDist]{g2}{g1}{$z$}
\boagate[mygreen][\DblVerDist]{g3}{g2}{$z$}
\broagate[myblue][\VerDist and 2*\HorDist]{g4}{g1}{$x$}
\broagate[myblue][\VerDist and \HorDist]{g5}{g2}{$x$}
\aroagate[mygreen][\VerDist and \HorDist]{g6}{g4}{$z$}
\aroagate[mygreen][\VerDist and 2*\HorDist]{g7}{g5}{$z$}
\boagate[mygreen][\DblVerDist]{g8}{g7}{$z$}
\broagate[myblue][\VerDist and 2*\HorDist]{g9}{g6}{$x$}
\broagate[myblue][\VerDist and \HorDist]{g10}{g7}{$x$}
\aroagate[mygreen][\VerDist and \HorDist]{g11}{g9}{$z$}
\aroagate[mygreen][\VerDist and 2*\HorDist]{g12}{g10}{$z$}
\boagate[mygreen][\DblVerDist]{g13}{g12}{$z$}
\broagate[myblue][\VerDist and 2*\HorDist]{g13}{g11}{$x$}
\broagate[myblue][\VerDist and \HorDist]{g14}{g12}{$x$}
\end{tikzpicture}
}\hfill%
\subfloat[Kitaev triangle]{%
\label{fig:TFK3}
\begin{tikzpicture}[on grid, auto]
\renewcommand{\ExtSize}{3.6mm} 
\renewcommand{\IntHeight}{9mm}
\renewcommand{\HorDist}{5.4mm}
\renewcommand{\VerDist}{5.4mm}
\renewcommand{\DblVerDist}{2*\VerDist} 
\renewcommand{\DblHorDist}{2*\HorDist}
\tikzstyle{oneaxis}=[draw,rectangle,thick,rounded corners = 1.2mm,minimum height=\IntHeight,adjust width=\IntAspect,inner sep=0,fill=white]
\fqwires[12*\HorDist]{w1}{6}
\aroagate[mygreen][\VerDist/2 and \HorDist]{g1}{w161}{$z$}
\aroagate[myblue]{g2}{g1}{$x$}
\aroagate[mygreen]{g3}{g2}{$z$}
\aroagate[myblue]{g4}{g3}{$x$}
\aroagate[mygreen][\VerDist and 2*\HorDist]{g5}{g4}{$x$}
\roagate[mygreen][\DblHorDist]{g6}{g1}{$z$}
\aroagate[myblue][\VerDist and 2*\HorDist]{g7}{g6}{$x$}
\aroagate[mygreen]{g8}{g7}{$z$}
\aroagate[myblue]{g9}{g8}{$x$}
\roagate[mygreen][3*\HorDist]{g10}{g6}{$z$}
\aroagate[myblue][\VerDist and 2*\HorDist]{g11}{g10}{$x$}
\aroagate[mygreen]{g12}{g11}{$z$}
\roagate[mygreen][3*\HorDist]{g13}{g10}{$z$}
\aroagate[myblue]{g14}{g13}{$x$}
\roagate[mygreen][\DblHorDist]{g15}{g13}{$z$}
\end{tikzpicture}
}\hfill~\\

~\hfill%
\subfloat[TFIM square]{%
\label{fig:TFK2}
\begin{tikzpicture}[on grid, auto]
\renewcommand{\ExtSize}{3.6mm} 
\renewcommand{\IntHeight}{9mm}
\renewcommand{\HorDist}{5.4mm}
\renewcommand{\VerDist}{5.4mm}
\renewcommand{\DblVerDist}{2*\VerDist} 
\renewcommand{\DblHorDist}{2*\HorDist}
\tikzstyle{oneaxis}=[draw,rectangle,thick,rounded corners = 1.2mm,minimum height=\IntHeight,adjust width=\IntAspect,inner sep=0,fill=white]
\fqwires[10*\HorDist]{w2}{3}
\rextgate[mygreen]{g21}{w211}{$z$}
\rextgate[mygreen]{g22}{w221}{$z$}
\rextgate[mygreen]{g23}{w231}{$z$}
\broagate[myblue][\VerDist/2 and \HorDist]{g24}{g22}{$x$}
\aroagate[myblue][\VerDist and \HorDist]{g25}{g24}{$x$}
\rextgate[mygreen][3*\HorDist]{g26}{g21}{$z$}
\bextgate[mygreen]{g27}{g26}{$z$}
\bextgate[mygreen]{g28}{g27}{$z$}
\broagate[myblue][\VerDist/2 and \HorDist]{g29}{g27}{$x$}
\aroagate[myblue][\VerDist and \HorDist]{g30}{g29}{$x$}
\rextgate[mygreen][3*\HorDist]{g31}{g26}{$z$}
\bextgate[mygreen]{g32}{g31}{$z$}
\bextgate[mygreen]{g33}{g32}{$z$}
\broagate[myblue][\VerDist/2 and \HorDist]{g34}{g32}{$x$}
\aroagate[myblue][\VerDist and \HorDist]{g35}{g34}{$x$}
\end{tikzpicture}
}\hfill%
\subfloat[TFIM triangle]{%
\label{fig:TFK4}
\begin{tikzpicture}[on grid, auto]
\renewcommand{\ExtSize}{3.6mm} 
\renewcommand{\IntHeight}{9mm}
\renewcommand{\HorDist}{5.4mm}
\renewcommand{\VerDist}{5.4mm}
\renewcommand{\DblVerDist}{2*\VerDist} 
\renewcommand{\DblHorDist}{2*\HorDist}
\tikzstyle{oneaxis}=[draw,rectangle,thick,rounded corners = 1.2mm,minimum height=\IntHeight,adjust width=\IntAspect,inner sep=0,fill=white]
\fqwires[12*\HorDist]{w2}{3}
\rextgate[mygreen]{g21}{w231}{$z$}
\aroagate[myblue][\VerDist/2 and \HorDist]{g22}{g21}{$x$}
\rextgate[mygreen][3*\HorDist]{g23}{w221}{$z$}
\aroagate[myblue][\VerDist/2 and \HorDist]{g24}{g23}{$x$}
\rextgate[mygreen][6*\HorDist]{g25}{w211}{$z$}
\rextgate[mygreen][2*\HorDist]{g26}{g21}{$z$}
\aroagate[myblue][\VerDist/2 and 2*\HorDist]{g27}{g26}{$x$}
\arextgate[mygreen][\VerDist/2 and \HorDist]{g28}{g27}{$z$}
\aroagate[myblue][\VerDist/2 and \HorDist]{g29}{g28}{$x$}
\rextgate[mygreen][3*\HorDist]{g30}{g26}{$z$}
\aroagate[myblue][\VerDist/2 and 2*\HorDist]{g31}{g30}{$x$}
\arextgate[mygreen][\VerDist/2 and \HorDist]{g32}{g31}{$z$}
\rextgate[mygreen][3*\HorDist]{g33}{g30}{$z$}
\aroagate[myblue][\VerDist/2 and \HorDist]{g34}{g33}{$x$}
\rextgate[mygreen][2*\HorDist]{g35}{g33}{$z$}
\end{tikzpicture}
}\hfill~\\
\label{fig:TFIMKitaev}
\caption{(a) A square circuit for a Kitaev chain with $2N = 6$ spins and (b)
the corresponding triangle circuit. These can be mapped to a (c) square circuit 
for a TFIM Hamiltonian for $N=3$ spins and (d) a triangle circuit for the TFIM Hamiltonian.}
\end{figure}

An alternative approach is to use \Cref{lem:TFIMto} as a turnover operation for a TFIM Hamiltonian.
In this setting, the turnover simultaneously operates on one- and two-spin gates and the 
compression algorithm appears to be different from the algorithms in \Cref{sec:algorithms},
but it turns out that they are completely analogous.
We illustrate the gist of the idea for a small example that draws a parallel between a $6$-spin
Kitaev chain and a $3$-spin TFIM circuit in \Cref{fig:TFIMKitaev}.
We see that we can map the two-spin $R^{zz}$ rotations from the Kitaev chain to the one-spin $R^{z}$
rotations in the TFIM Hamiltonian and the $R^{xx}$ rotations are mapped to $R^{xx}$ rotations 
that mutually commute (\Cref{lem:pcom}).
Another interpretation is that every pair of consecutive spins in the Kitaev chain is combined to a 
single spin in the TFIM circuit.
It follows from this mapping that we can use the algorithms from \Cref{sec:algorithms}
together with \Cref{lem:pfus,lem:TFIMto} to compress TFIM circuits and transform
them from square to triangle or vice versa.
The total number of gates for a minimal representation of an $N$-spin TFIM circuit is $N (2N-1)$,
$N^2$ of the gates are one-spin Pauli rotations, $N(N-1)$ are two-spin Pauli rotations.

Depending on the details of the quantum hardware, it is important to remark that the cost
of a single two-spin Pauli-X rotation can sometimes be considered approximately the same as the cost of a two-spin
XY rotation since both require the same number of two qubit CNOT gates~\cite{Bassman2021}.
These are the main source of errors on many devices and in that case we expect the TFXY mapping to perform
better.
However, for other devices the $R^{\alpha\alpha}$ transformations can be supported natively and the TFIM
mapping might be preferable.


\section{Implementation details}
\label{sec:implementation}

We describe some relevant details of the numerical implementations
of our fusion and turnover algorithms.

\subsection{Storage}

All Pauli rotations throughout the compression algorithms for 
Ising, Kitaev chains, XY, and TFIM models are stored by two doubles
storing $\cos(\theta/2)$ and $\sin(\theta/2)$.
This avoids the numerical evaluation of (inverse) trigonometric functions,
and instead we can rely on Givens rotation matrices in the turnover operations.

For the fusion of two compatible 1- or 2-spin Pauli rotations (\Cref{lem:pfus}),
we simply compute the first column of a product of two $2 \times 2$ matrices.

\subsection{Turnover of SU(2)}

The angles for an $\SU2$ turnover are given by \eqref{eq:toSU21,eq:toSU22}.
Numerical evaluation of these formulas requires (inverse) trigonometric functions
which can be avoided and replaced by Givens rotations leading to ideal 
numerical roundoff properties \cite{Aurentz2018}.

As input to our $\SU2$ turnover routine, we have three Euler angles stored as
\begin{align}
& \begin{bmatrix} \cos(\nicefrac{\theta_a}{2}) & \sin(\nicefrac{\theta_a}{2}) \end{bmatrix}, &
& \begin{bmatrix} \cos(\nicefrac{\theta_b}{2}) & \sin(\nicefrac{\theta_b}{2}) \end{bmatrix}, &
\begin{bmatrix} \cos(\nicefrac{\theta_c}{2}) & \sin(\nicefrac{\theta_c}{2}) \end{bmatrix},
\label{eq:inputEuler}
\end{align}
or, in short, $[c_a, s_a]$, $[c_b, s_b]$, and $[c_c, s_c]$.
Without loss of generality, we form the first column of the corresponding $\SU2$ matrix \eqref{eq:SU2}
under the assumption that the Euler angles correspond to a YZY parametrization \eqref{eq:yzy}, $R^y(\theta_1) R^z(\theta_2) R^y(\theta_3)$.
The $(1,1)$ element, $\alpha = \alpha_r + \I \alpha_i$, and the $(2,1)$ element, $\beta = \beta_r + \I \beta_i$, of the $\SU2$ matrix are in that case:
\begin{equation}
\begin{aligned}
\alpha_r & = c_b ( c_a c_c - s_a s_c ), & \alpha_i & = - s_b ( c_a c_c + s_a s_c ), \\
\beta_r & = c_b ( s_a c_c + c_a s_c ), & \beta_i & = - s_b ( s_a c_c - c_a s_c ), 
\label{eq:yzynum}
\end{aligned}
\end{equation}
which we directly compute from our input Euler angles \eqref{eq:inputEuler}.
The dual Euler angles that we have to compute form a ZYZ decomposition \eqref{eq:zyz} of the same $\SU2$
matrix, so they have to sayisfy:
\begin{equation}
\begin{aligned}
\alpha_r & = c_2 ( c_1 c_3 - s_1 s_3 ), & \alpha_i & = - c_2 ( s_1 c_3 + c_1 s_3 ), \\
\beta_r & = s_2 ( c_1 c_3 + s_1 s_3 ), & \beta_i & = \phantom{-}s_2 ( s_1 c_3 - c_1 s_3 ).
\label{eq:zyznum}
\end{aligned}
\end{equation}
It follows from \eqref{eq:zyz,eq:zyznum} that we can compute $c_2 = |\alpha| = \sqrt{ \alpha_r^2 + \alpha_i^2}$ and $s_2 = |\beta| = \sqrt{ \beta_r^2 + \beta_i^2 }$, which is properly normalized as
$|\alpha|^2 + |\beta|^2 = 1$ by \eqref{eq:SU2}.

Next, we get from \eqref{eq:zyznum} that we can compute $[c_1, \, s_1]$ and $[c_3, \, s_3]$
from two Givens rotation
matrices that introduce zeros in, respectively, $v_1$ or $v_1^{\prime}$, and $v_3$ or $v_3^{\prime}$,
\begin{equation}
\begin{aligned}
v_1 & = 
\begin{bmatrix}
\alpha_r / c_2 + \beta_r / s_2 \\
-\alpha_i / c_2 + \beta_i/ s_2
\end{bmatrix}, &
v_1^{\prime} & =
\begin{bmatrix}
-\alpha_i / c_2 - \beta_i / s_2 \\
\beta_r / s_2 - \alpha_r / c_2
\end{bmatrix}, \\
v_3 & = 
\begin{bmatrix}
-\alpha_i / c_2 + \beta_i /s_2 \\
-\alpha_r / c_2 + \beta_r/s_2
\end{bmatrix}, &
v_3^{\prime} & = 
\begin{bmatrix}
\alpha_r / c_2 + \beta_r /s_2 \\
-\alpha_i / c_2 - \beta_i/s_2
\end{bmatrix}.  
\end{aligned}
\end{equation}
To ensure the lowest relative error, we choose the vector with largest
norm between $v_1$ and $v_1^{\prime}$, and similar for $v_3$ and $v_3^{\prime}$.
It follows from the roundoff properties of Givens rotation matrices \cite{Aurentz2018}
that this algorithm computes the dual Euler angles $[ c_1, s_1 ]$, $[ c_2, s_2 ]$ 
and $[ c_3, s_3 ]$ to high relative accuracy. 

\subsection{TFXY gates}

A TFXY gate is parametrized by six angles and has the following
matrix representation,
\begin{equation}
U_{\text{TFXY}} =
\begin{tikzpicture}[on grid, auto, baseline=-\VerDist/2]
\fqwires[4*\HorDist]{w1}{2}
\rextgate[mygreen]{gate2}{w111}{$z$}
\rextgate[mygreen]{gate3}{w121}{$z$}
\brtagate[\VerDist/2 and 2*\HorDist]{myblue,myred}{gate1}{w111}{$x$}{$y$}
\lextgate[mygreen]{gate4}{w112}{$z$}
\lextgate[mygreen]{gate5}{w122}{$z$}
\anode[3*\HorDist/4]{p1}{gate2}{\footnotesize$\theta_1$}
\bnode[3*\HorDist/4]{p2}{gate3}{\footnotesize$\theta_2$}
\anode[5*\HorDist/4]{p1}{gate1}{\footnotesize$\theta_3$}
\bnode[5*\HorDist/4]{p2}{gate1}{\footnotesize$\theta_4$}
\anode[3*\HorDist/4]{p1}{gate4}{\footnotesize$\theta_5$}
\bnode[3*\HorDist/4]{p2}{gate5}{\footnotesize$\theta_6$}
\end{tikzpicture}
 = 
\begin{bmatrix}
\alpha & & & - \bar{\delta} \\
 & \beta & -\bar{\gamma} & \\
 & \gamma & \bar{\beta} & \\
\delta & & & \bar{\alpha} 
\end{bmatrix},
\label{eq:TFXYangles}
\end{equation}
where the matrix entries are given by
\begin{equation}
\begin{aligned}
\alpha & = \phantom{-\I} \cos((\theta_3 - \theta_4)/2) \ e^{-\I(\theta_1 + \theta_2 + \theta_5 + \theta_6)/2},\\
\beta & = \phantom{-\I}  \cos((\theta_3 + \theta_4)/2) \ e^{-\I(\theta_1 - \theta_2 + \theta_5 - \theta_6)/2},\\
\gamma & = -\I \sin((\theta_3 + \theta_4)/2) \ e^{-\I(\theta_1 - \theta_2 - \theta_5 + \theta_6)/2},\\
\delta & = -\I \sin((\theta_3 - \theta_4)/2) \ e^{-\I(\theta_1 + \theta_2 - \theta_5 - \theta_6)/2}.
\end{aligned}
\label{eq:TFXYelem}
\end{equation} 
It can be shown that both the outer $2 \times 2$ matrix
$\left[ \begin{smallmatrix} \alpha & -\bar{\delta}\\ \delta & \bar{\alpha} \end{smallmatrix}\right]$
and the inner $2 \times 2$ matrix
$\left[ \begin{smallmatrix} \beta & -\bar{\gamma}\\ \gamma & \bar{\beta} \end{smallmatrix}\right]$
form $\SU2$ by mapping the matrix elements to two independent ZYZ decompositions \eqref{eq:zyz}.
In our implementation, we store a TFXY gate by the four complex numbers $\alpha$, $\beta$, 
$\gamma$, and $\delta$ which can be easily computed form the angles via \eqref{eq:TFXYelem}.
The values $\alpha$, $\beta$, $\gamma$, and $\delta$ 
can easily be converted back to a parametrization with six angles.

Unitary $4 \times 4$ matrices that have the nonzero pattern of the matrix in \eqref{eq:TFXYangles}
are also known as \emph{matchgates}\cite{Bassman2021}. Matchgates can be permuted to 
a block-diagonal matrix by the following permutation matrix:
\begin{equation}
P_{\oplus} = 
\begin{bmatrix}
1 & & & \\ & & & 1 \\ & 1 & &  \\ &  & 1 & 
\end{bmatrix},
\qquad
P_{\oplus}^{T} =
\begin{bmatrix}
1 & & & \\ & & 1 & \\ & & & 1 \\ & 1 & & 
\end{bmatrix}.
\label{eq:permutation}
\end{equation}
We get that, 
\begin{equation}
P_{\oplus}
\begin{bmatrix}
\alpha & & & - \bar{\delta} \\
 & \beta & -\bar{\gamma} & \\
 & \gamma & \bar{\beta} & \\
\delta & & & \bar{\alpha} 
\end{bmatrix}
P_{\oplus}^{T} =
\begin{bmatrix}
\alpha & -\bar{\delta} & & \\
\delta & \bar{\alpha} & & \\
 & & \beta & -\bar{\gamma} \\
 & & \gamma & \bar{\beta}
\end{bmatrix}
=
\begin{bmatrix} \alpha & -\bar{\delta}\\ \delta & \bar{\alpha} \end{bmatrix}
\oplus
\begin{bmatrix} \beta & -\bar{\gamma}\\ \gamma & \bar{\beta} \end{bmatrix}.
\label{eq:mtod}
\end{equation}
In what follows, we use the notation $\TFXY(A,B)$, where $A, B \in \SU2$, to denote
the TFXY gate with its outer $2 \times 2$ $\SU2$ matrix equal to $A$ and its middle $2 \times 2$
$\SU2$ matrix equal to $B$.
With this notation, we can rewrite \eqref{eq:mtod} in schematic notation as
\begin{equation*}
{\begin{myqcircuit*}{1}{0.75}
  & \multigate1{P_{\oplus}^T} & \multigate1{\TFXY(A,B)} & \multigate1{P_{\oplus}} & \qw\\
  & \ghost{P_{\oplus}^T}      & \ghost{\TFXY(A,B)}      & \ghost{P_{\oplus}}      & \qw
\end{myqcircuit*}}
\ = \
{\begin{myqcircuit*}{1}{0.75}
  & \multigate1{A \oplus B} & \qw\\
  & \ghost{A \oplus B}      & \qw
\end{myqcircuit*}} \quad .
\end{equation*}

The fusion of two TFXY gates can thus simply be computed 
as the matrix product of two $2 \times 2$ matrices:
\begin{equation*}
\TFXY(A,B) = \TFXY(A_1,B_1) \, \TFXY(A_2,B_2),
\quad \Longleftrightarrow \quad A = A_1 A_2, \ B = B_1 B_2.
\end{equation*}
%
Two sets of three TFXY gates satisfy a turnover relation if,
\begin{equation}
\resizebox{.9 \textwidth}{!} 
{$
{\begin{myqcircuit*}{1}{0.75}
  & \multigate1{\TFXY(A,B)} & \qw                    & \multigate1{\TFXY(E,F)} & \qw\\
  & \ghost{\TFXY(A,B)}      & \multigate1{\TFXY(C,D)} & \ghost{\TFXY(E,F)}.    & \qw\\
  & \qw                    & \ghost{\TFXY(C,D)}      & \qw                     & \qw
\end{myqcircuit*}}
\ = \
{\begin{myqcircuit*}{1}{0.75}
  & \qw                     & \multigate1{\TFXY(W,X)} & \qw                     & \qw \\
  & \multigate1{\TFXY(U,V)} & \ghost{\TFXY(W,X)}      & \multigate1{\TFXY(Y,Z)} & \qw\\
  & \ghost{\TFXY(U,V)}      & \qw                     & \ghost{\TFXY(Y,Z)}      & \qw    
\end{myqcircuit*}} \quad .
$}
\label{eq:TFXYTO}
\end{equation}
We will present how we numerically compute the turnover transformation for
TFXY matrices.
For every $A, B \in \SU2$, we have that,
\begin{equation}
{\begin{myqcircuit*}{1}{0.75}
  & \qw                     & \qw                     & \qw                       & \qw\\
  & \multigate1{P_{\oplus}^T} & \multigate1{\TFXY(A,B)} & \multigate1{P_{\oplus}} & \qw\\
  & \ghost{P_{\oplus}^T}      & \ghost{\TFXY(A,B)}      & \ghost{P_{\oplus}}      & \qw
\end{myqcircuit*}}
\ = \
{\begin{myqcircuit*}{1}{0.75}
  & \qw                     & \qw\\
  & \multigate1{A \oplus B} & \qw\\
  & \ghost{A \oplus B}      & \qw
\end{myqcircuit*}}
\ = \ A \oplus B \oplus A \oplus B,
\label{eq:perm1}
\end{equation}
which is a block-diagonal matrix with repeated blocks on the diagonal.
Similarly, from carrying out the matrix multiplication, we find that,
\begin{equation}
{\begin{myqcircuit*}{1}{0.75}
  & \qw                     & \multigate1{\TFXY(A,B)} & \qw                       & \qw\\
  & \multigate1{P_{\oplus}^T} & \ghost{\TFXY(A,B)}    & \multigate1{P_{\oplus}} & \qw\\
  & \ghost{P_{\oplus}^T}      & \qw                   & \ghost{P_{\oplus}}      & \qw
\end{myqcircuit*}}
\ = \
\left[
\begin{smallmatrix}
\alpha & & & & & & & -\bar{\delta} \\
& \beta & & & & & -\bar{\gamma} & \\
& & \alpha & & & -\bar{\delta} & & \\
& & & \beta & -\bar{\gamma} & & & \\
& & & \gamma & \bar{\beta} & & & \\
& & \delta & & & \bar{\alpha} & & \\
& \gamma & & & & & \bar{\beta} & \\
\delta & & & & & & & \bar{\alpha} \\
\end{smallmatrix}\right],
\label{eq:perm2}
\end{equation}
for $A = \left[ \begin{smallmatrix} \alpha & -\bar{\delta}\\ \delta & \bar{\alpha} \end{smallmatrix}\right]$
and $B = \left[ \begin{smallmatrix} \beta & -\bar{\gamma}\\ \gamma & \bar{\beta} \end{smallmatrix}\right]$.

It follows from \eqref{eq:perm1} and \eqref{eq:perm2} that in the permuted 
basis, $\eye[2] \otimes P_{\oplus}$, the TFXY turnover equation \eqref{eq:TFXYTO}
becomes a matrix equation involving matrices of the following nonzero structure:
\begin{equation*}
\resizebox{\textwidth}{!}{%
\begin{tikzpicture}
\crossmat{0}{0}
\anode[40]{AB}{m44}{$A,B$}
\blockdiagmat{24}{0}
\anode[40]{AB}{m44}{$C,D$}
\crossmat{48}{0}
\anode[40]{AB}{m44}{$E,F$}
\rnode[48]{eq1}{m55}{$=$}
\blockdiagmat{80}{0}
\anode[40]{AB}{m44}{$U,V$}
\crossmat{104}{0}
\anode[40]{AB}{m44}{$W,X$}
\blockdiagmat{128}{0}
\anode[40]{AB}{m44}{$Y,Z$}
\end{tikzpicture}
}
\end{equation*}
The squares and circles indicate the elements from each $\SU2$ matrix
that make up the TFXY gate.
To go from the factorization on the left-hand side to the factorization on the
right-hand side ($\VEE$ to $\HAT$), we form the product of the three unitaries on the left-hand side
and block-divide this unitary matrix $Q$ as
\begin{equation*}
\begin{bmatrix}
Q_{11} & & & Q_{14}\\
 & Q_{22} & Q_{23} & \\
 & Q_{32} & Q_{33} & \\
Q_{41} & & & Q_{44} 
\end{bmatrix}
= \quad
\scalebox{0.7}{%
\begin{tikzpicture}[baseline=-11mm]
\turnovermat{0}{0}{black}
\end{tikzpicture}}.
\end{equation*}
This unitary matrix has a lot of structure that we exploit.
The four $2 \times 2$ matrix pairs $(Q_{11}, Q_{33})$, $(Q_{22}, Q_{44})$, 
$(Q_{23}, Q_{41})$, $(Q_{14}, Q_{32})$ are all of the form:
\begin{equation}
\left(
\begin{bmatrix}
a & b \\
c & d
\end{bmatrix},
\begin{bmatrix}
\bar{d} & -\bar{c} \\
-\bar{b} & \bar{a}
\end{bmatrix}
\right).
\end{equation}
This means that we only have to compute four $2 \times 2$ matrices
to construct $Q$.
Furthermore, it can be shown that matrix pairs of this form can always
be simultaneously (anti-)diagonalized.
A second property that we use, is that the four diagonal blocks are of equal norm, and similarly the four blocks on the anti-diagonal are of equal norm.
To compute the turnover, we have to compute the following four coupled matrix decompositions:
\begin{equation*}
\begin{aligned}
U^* (Q_{11}, Q_{33}) Y^* & = \left(\left[\begin{smallmatrix} w_{11} & \\ & x_{11} \end{smallmatrix}\right], \left[\begin{smallmatrix} x_{22} & \\ & w_{22} \end{smallmatrix}\right]\right), &
V^* (Q_{22}, Q_{44}) Z^* & = \left(\left[\begin{smallmatrix} w_{11} & \\ & x_{11} \end{smallmatrix}\right], \left[\begin{smallmatrix} x_{22} & \\ & w_{22} \end{smallmatrix}\right]\right), \\
U^* (Q_{14}, Q_{32}) Z^* & = \left( \left[\begin{smallmatrix} & w_{12} \\ x_{12} & \end{smallmatrix}\right], \left[\begin{smallmatrix} & x_{21} \\ w_{21} & \end{smallmatrix}\right]\right), &
V^* (Q_{23}, Q_{41}) Y^* & = \left( \left[\begin{smallmatrix} & w_{12} \\ x_{12} & \end{smallmatrix}\right], \left[\begin{smallmatrix} & x_{21} \\ w_{21} & \end{smallmatrix}\right]\right),
\end{aligned}
\label{eq:diagonalization}
\end{equation*}
in such a way that the ordering of the (anti-)diagonalizations is consistent and that they share
the same eigenvectors.

Our numerical algorithm depends on the ratio between $\|Q_{i,i}\|_{2}$ and $\|Q_{k,5-k}\|_{2}$.
If $\|Q_{i,i}\|_{2} \leq \|Q_{k,5-k}\|_{2}$, we start with computing $\SU2$ matrices
$U$ and $Y$ that diagonalize $(Q_{11}, Q_{33})$. Keeping $Y$ fixed, we compute $V$
that anti-diagonalizes $(Q_{23}, Q_{41})Y^*$. Finally keeping $V$ fixed, we compute $Z$
that diagonalizes $V^*(Q_{22},Q_{44})$.
If $\|Q_{i,i}\|_{2} > \|Q_{k,5-k}\|_{2}$, we start with computing $\SU2$ matrices $V$ and $Y$ that
anti-diagonalize $(Q_{23}, Q_{41})$. Keeping $Y$ fixed, we compute $U$ that diagonalizes
$(Q_{11},Q_{33})Y^*$. Finally keeping $U$ fixed, we compute $Z$ that anti-diagonalizes
$U^*(Q_{14},Q_{32})$.

Afterwards, the $\SU2$ matrices $W$ and $X$ are determined from the (anti-)diagonal elements 
of the (anti-)diagonalized matrices.
To compute the turnover in the other direction ($\HAT$ to $\VEE$), 
from the right-hand side of \eqref{eq:TFXYTO} to the left-hand side, 
we simply interchange the roles of the first and third spins.
This is just a permutation of the $8 \times 8$ unitary.
During our extensive numerical tests we observed that this algorithm always computes
the TFXY turnover up to high relative accuracy.


\section{Numerical examples}
\label{sec:ex}

All numerical experiments are performed on a
AMD Ryzen Threadripper 3990X 64-Core Processor @ 2.9 GHz with 256 GB RAM.
Our experiments can be reproduced with the \texttt{F3C++} code \cite{f3c,f3cpp}.

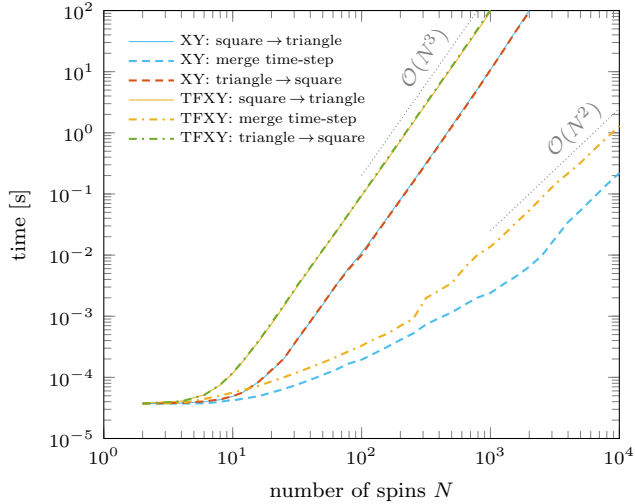
\begin{figure}[hbtp]
\centering
\figname{scaling}%
\begin{tikzpicture}
\begin{loglogaxis}[
  xmin=1,xmax=1e4,%
  ymin=1e-5,ymax=1e2,%
  xlabel={number of spins $N$},%
  ylabel={time [s]},%
  legend style={draw=none,fill=none,row sep=-2pt},%
  legend pos=north west,%
]
\addplot[myColSix,mark=]%
   table[x index=0,y index=1] {\datfile{XY-square2triangle}};
\addplot[myColSix,thick,densely dashed,mark=]%
   table[x index=0,y index=1] {\datfile{XY-merge-timestep}};
\addplot[myColTwo,thick,densely dashed,mark=]%
   table[x index=0,y index=1] {\datfile{XY-triangle2square}};
\addplot[myColThr,mark=]%
   table[x index=0,y index=1] {\datfile{TFXY-square2triangle}};
\addplot[myColThr,thick,dash dot,mark=]%
   table[x index=0,y index=1] {\datfile{TFXY-merge-timestep}};
\addplot[myColFiv,thick,dash dot,mark=]%
   table[x index=0,y index=1] {\datfile{TFXY-triangle2square}};
\legend{\scalebox{0.8}{XY: square\,$\to$\,triangle},
        \scalebox{0.8}{XY: merge time-step},
        \scalebox{0.8}{XY: triangle\,$\to$\,square},
        \scalebox{0.8}{TFXY: square\,$\to$\,triangle},
        \scalebox{0.8}{TFXY: merge time-step},
        \scalebox{0.8}{TFXY: triangle\,$\to$\,square}};
\addplot[no marks,gray,densely dotted] plot coordinates {(1e2,2e-1) (1e3,2e2)};
\addplot[no marks,gray,densely dotted] plot coordinates {(1e3,2.5e-2) (1e4,2.5e0)};
\draw[gray] (3.5e2,1e2) node[below left,rotate=55] {\footnotesize$\bigO(N^3)$};
\draw[gray] (1e4,2e0) node[above left,rotate=44] {\footnotesize$\bigO(N^2)$};
\end{loglogaxis}
\end{tikzpicture}%
\vspace{-5pt}
\caption{\label{fig:scaling}%
Scaling for compression of XY and TFXY models in function of number of spins.}
\end{figure}

\begin{figure}[hbtp]
\centering
\begin{tikzpicture}
\begin{loglogaxis}[
  name=ax1,%
  width=.48\textwidth,%
  xmin=1,xmax=1e4,%
  ymin=1e-16,ymax=1e0,%
  ytick={1e-16,1e-12,1e-8,1e-4,1e0},%
  extra y ticks={1e-15,1e-14,1e-13,1e-11,1e-10,1e-9,1e-7,1e-6,1e-5,1e-3,1e-2,1e-1},%
  extra y tick labels={},%
  extra y tick style={tickwidth=\pgfkeysvalueof{/pgfplots/minor tick length}},%
  xlabel={number of time-steps $n_t$},%
  ylabel={$\|Q_{c} - Q_{f}\|_F$},%
  ylabel shift=-4pt,%
  legend style={draw=none,fill=none,row sep=-2pt},%
  legend pos=north west,%
  legend cell align=right,%
]
\draw (1e4,1e0) node[below left,xshift=-5pt,yshift=-6pt] {XY};
\addlegendentry{6 spins}
\addplot[myColOne,mark=]%
   table[x index=0,y index=1] {\datfile{XY-error-6q}};
\addlegendentry{8 spins}
\addplot[myColTwo,mark=]%
   table[x index=0,y index=1] {\datfile{XY-error-8q}};
\addlegendentry{10 spins}
\addplot[myColThr,mark=]%
   table[x index=0,y index=1] {\datfile{XY-error-10q}};
\addplot[no marks,gray,densely dotted] plot coordinates {(1e0,1e-12) (1e4,1e-8)};
\draw[gray] (1e4,1e-8) node[above left,rotate=12] {\footnotesize$\bigO(n_t)$};
\end{loglogaxis}
%
\begin{loglogaxis}[
  at={(ax1.south east)},%
  xshift=2cm,%
  width=.48\textwidth,%
  xmin=1,xmax=1e4,%
  ymin=1e-16,ymax=1e0,%
  ytick={1e-16,1e-12,1e-8,1e-4,1e0},%
  extra y ticks={1e-15,1e-14,1e-13,1e-11,1e-10,1e-9,1e-7,1e-6,1e-5,1e-3,1e-2,1e-1},%
  extra y tick labels={},%
  extra y tick style={tickwidth=\pgfkeysvalueof{/pgfplots/minor tick length}},%
  xlabel={number of time-steps $n_t$},%
  ylabel={$\|Q_{c} - Q_{f}\|_F$},%
  ylabel shift=-4pt,%
  legend style={draw=none,fill=none,row sep=-2pt},%
  legend pos=north west,%
  legend cell align=right,%
]
\draw (1e4,1e0) node[below left,xshift=-5pt,yshift=-6pt] {TFXY};
\addlegendentry{6 spins}
\addplot[myColOne,mark=]%
   table[x index=0,y index=1] {\datfile{TFXY-error-6q}};
\addlegendentry{8 spins}
\addplot[myColTwo,mark=]%
   table[x index=0,y index=1] {\datfile{TFXY-error-8q}};
\addlegendentry{10 spins}
\addplot[myColThr,mark=]%
   table[x index=0,y index=1] {\datfile{TFXY-error-10q}};
\addplot[no marks,gray,densely dotted] plot coordinates {(1e0,1e-12) (1e4,1e-8)};
\draw[gray] (1e4,1e-8) node[above left,rotate=12] {\footnotesize$\bigO(n_t)$};
\end{loglogaxis}
\end{tikzpicture}%
\vspace{-5pt}%
\caption{\label{fig:error}%
Maximum Frobenius error on 100 randomly generated compressed circuits ($Q_c$)
compared to full circuits for XY and TFXY models
in function of number of time-steps.}
\end{figure}
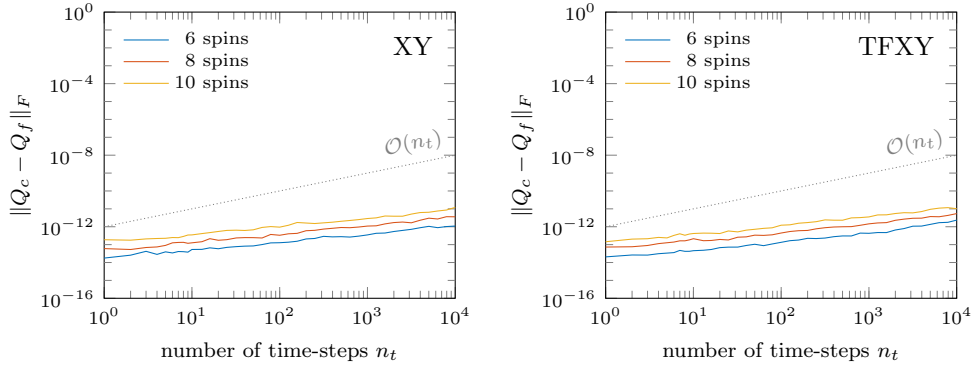

For our first numerical experiment, we compressed Trotter circuits
for XY and TFXY models with randomly generated angles sampled from the standard normal distribution
to benchmark the speed and accuracy of our algorithms.
\Cref{fig:scaling} shows the cubic scaling of transforming the circuit between
square and triangle representation and back, and the quadratic scaling of merging
a time-step with a triangle circuit.
Both algorithms easily scale up to $\bigO(10^3)$ spins, which is well beyond the capabilities
of current quantum hardware.
The compression of a full TFXY model, based on the simultaneous diagonalization 
algorithm for the turnover operation, is about
a factor of 10 slower than compressing an XY model, which uses two $\SU2$ turnovers.
The TFXY turnover can alternatively be implemented with 32 $\SU2$ turnovers \cite{TrotterCompression},
hence the simultaneous diagonalization approach is 1.5 times faster.
\Cref{fig:error} displays the Frobenius norm difference between the unitary of the full Trotter
circuit $Q_f$, and the unitary of the compressed Trotter circuit $Q_c$ in function of 
the number of time-steps for systems up to 10 spins.
We observe that our compression algorithm achieves high accuracy in all cases and that the
rounoff error grows sub-linearly.

For our final numerical experiment, we report the error on the compression algorithms
for the adiabatic state preparation experiment presented in~\cite{TrotterCompression}.
In this experiment, the ground state of a classical Ising model, that is easy to prepare,
is evolved adiabatically to the ground state of a TFIM Hamiltonian.
The experiment in~\cite{TrotterCompression} is performed for a 5-spin system with two
different time-steps $\Delta t$, $0.05$ and $0.25$, and $n_t$ respectively chosen as $1200$ and
$240$ such that the final simulation time is the same. The first option leads to a significantly
smaller Trotter error but to a slightly higher numerical error on the compressed circuit due to roundoff
accumulation as five times more time-steps have to be merged.
The total error remains small in all cases.

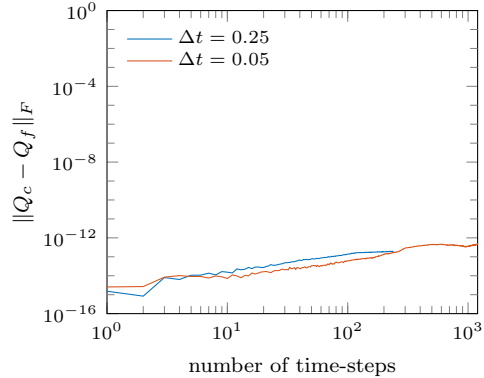
\begin{figure}[hbtp]
\centering
\begin{tikzpicture}
\begin{loglogaxis}[
  width=0.5\textwidth,
  xmin=1,%
  xmax=1200,%
  ymin=1e-16,ymax=1e0,%
  ytick={1e-16,1e-12,1e-8,1e-4,1e0},%
  extra y ticks={1e-15,1e-14,1e-13,1e-11,1e-10,1e-9,1e-7,1e-6,1e-5,1e-3,1e-2,1e-1},%
  extra y tick labels={},%
  extra y tick style={tickwidth=\pgfkeysvalueof{/pgfplots/minor tick length}},%
  xlabel={number of time-steps},%
  ylabel={$\|Q_{c} - Q_{f}\|_F$},%
  ylabel shift=-4pt,%
  legend style={draw=none,fill=none,row sep=-2pt},%
  legend pos=north west,%
]
\addlegendentry{$\Delta t = 0.25$}
\addplot[myColOne,mark=]%
   table[x index=0,y index=1] {\datfile{ASP-240}};
\addlegendentry{$\Delta t = 0.05$}
\addplot[myColTwo,mark=]%
   table[x index=0,y index=1] {\datfile{ASP-1200}};
\end{loglogaxis}
\end{tikzpicture}%
\vspace{-5pt}
\caption{\label{fig:ASP}%
Frobenius error on compressed circuit ($Q_c$) compared to full circuits ($Q_f$) for the
adiabatic state preparation experiment presented in~\cite{TrotterCompression}.}
\end{figure}


\section{Conclusions}
\label{sec:concl}

We presented a fast and accurate quantum circuit synthesis
algorithm that is suitable for compressing circuits
to simulate certain classes of spin Hamiltonians known as
free fermionizable models on current generation quantum devices.
Our algorithms are based on matrix factorizations and easily scale up 
to $\bigO(10^3)$ spins as we make use of localized circuit transformations
that allow us to keep the $2^N \times 2^N$ unitary matrix 
in factorized form throughout the algorithm.
Furthermore, we can make the Trotter error
as small as required without increasing the circuit depth.
Numerical experiments showed that our turnover routines based
on Givens rotation matrices and simultaneous (anti-)diagonalization
behave backward stable.
Our approach is not limited to Hamiltonian simulation, any quantum
circuit partially comprised of the gates that we discussed
admits local fusion and turnover operations.
We provide reference implementations that are readily available and
can be used to improve existing quantum compilers and transpilers.


\bibliographystyle{siamplain}
\bibliography{references}

\end{document}